%% file: A_Gour_Sahi_main_text8.tex
\documentclass[12pt]{article}
\usepackage{amsmath,amsthm}
\usepackage{amssymb,latexsym}

\title{On an analytic description of the $\alp$-cosine transform on real Grassmannians.}
\date{\today}
\author{Semyon Alesker\footnote{Partially supported by ISF grant 1447/12.}\\
{ \normalsize Department of Mathematics}\\
{ \normalsize Tel Aviv University, Ramat Aviv}\\
{ \normalsize 69978 Tel Aviv, Israel}\\
{ \normalsize e-mail: semyon@post.tau.ac.il}
\and
 Dmitry Gourevitch\\
{ \normalsize  Faculty of Mathematics and Computer Science}\\
{ \normalsize  Weizmann Institute of Science}\\
{ \normalsize  POB 26, Rehovot 76100, Israel}\\
{ \normalsize  e-mail: dimagur@weizmann.ac.il}
 \and
  Siddhartha Sahi\\
{ \normalsize   Department of Mathematics}\\
{ \normalsize   Rutgers University, Hill Center -Busch Campus}\\
{ \normalsize   110 Frelinghuysen Road Piscataway, NJ 08854-8019, USA}\\
{ \normalsize   e-mail: sahi@math.rugers.edu}}

\def\RR{\mathbb{R}}
\def\CC{\mathbb{C}}

\def\NN{\mathbb{N}}
\def\ZZ{\mathbb{Z}}

\def\PP{\mathbb{P}}

\def\FF{\mathbb{F}}

\def\inj{\hookrightarrow}

\def\eps{\varepsilon}
\def\alp{\alpha}

\def\lam{\lambda}

\def\to{\longrightarrow}

\def\qed { Q.E.D. }

\swapnumbers
\newtheorem{theorem}{Theorem}[section]
\newtheorem{corollary}[theorem]{Corollary}
\newtheorem{lemma}[theorem]{Lemma}
\newtheorem{proposition}[theorem]{Proposition}
\newtheorem{claim}[theorem]{Claim}
\theoremstyle{definition}
\newtheorem{example}[theorem]{Example}
\newtheorem{definition}[theorem]{Definition}
\newtheorem{remark}[theorem]{Remark}
 \def\cb{{\cal B}} 
\def\cd{{\cal D}}  
 \def\ch{{\cal H}} 
  \def\cl{{\cal L}}
 \def\cn{{\cal N}} \def\co{{\cal O}}

  \def\cu{{\cal U}}
\def\cv{{\cal V}}  
 
\def\fa{\mathfrak{a}}\def\fb{\mathfrak{b}}

\def\fg{\mathfrak{g}}\def\fh{\mathfrak{h}}
\def\fk{\mathfrak{k}}
\def\fm{\mathfrak{m}}\def\fn{\mathfrak{n}}
\def\fp{\mathfrak{p}}\def\fq{\mathfrak{q}}

\input diagram

\def\pt{\partial}

\numberwithin{equation}{section}
\begin{document}

\maketitle
\def\gk{G_k(V)}

\begin{abstract}
The goal of this paper is to describe the $\alp$-cosine transform on functions on real Grassmannian $Gr_i(\RR^n)$
in analytic terms as explicitly as possible. We show that for all but finitely many complex $\alp$ the $\alp$-cosine transform
is a composition of the $(\alp+2)$-cosine transform with an explicitly written (though complicated) $O(n)$-invariant differential operator.
For all exceptional values of $\alp$ except one we interpret the $\alp$-cosine transform explicitly as either the Radon transform
or composition of two Radon transforms. Explicit interpretation of the transform corresponding to the last remaining value $\alp$, which is $-(\min\{i,n-i\}+1)$, is still an open problem.
\end{abstract}

\tableofcontents

\section{Introduction}\label{S:Introduction}
The goal of this paper is to describe the $\alp$-cosine transform on functions on real Grassmannians
in analytic terms as explicitly as possible. This goal has been achieved in this paper to a certain extent, though
still there are questions requiring further clarification.

To formulate the problem more precisely, we have to remind the definition of the $\alpha$-cosine transform.
Fix an $n$-dimensional Euclidean space $V$. Let $1\leq i\leq n-1$ be an integer. Let $Gr_i(V)$ denote the Grassmannian
of $i$-dimensional linear subspaces. For a pair of subspaces $E,F\in Gr_i(V)$ one defines (the absolute value of)
the cosine of the angle $|\cos(E,F)|$ between $E$ and $F$ as the coefficient of the distortion of measure under
the orthogonal projection from $E$ to $F$. More precisely, let $q\colon E\to F$ denote the restriction to $E$ of
the orthogonal projection $V\to F$. Let $A\subset E$ be an arbitrary subset of finite positive Lebesgue measure.
Then define
$$|\cos(E,F)|:=\frac{vol_F(q(A))}{vol_E(A)},$$
where $vol_F, vol_E$ are Lebesgue measures on $F,E$ respectively induced by the Euclidean metric on $V$ normalized
so that the Lebesgue measures of unit cubes are equal to 1. It is easy to see that $|\cos(E,F)|$
is independent of the set $A$ and is symmetric with respect to $E$ and $F$.

\begin{example}
If $i=1$ then $|\cos(E,F)|$ is the usual cosine of the angle between two lines.
\end{example}

For $\alp\in \CC$ with $Re(\alp)\geq 0$ let us define the $\alpha$-cosine transform
$$T_\alp\colon C^\infty(Gr_i(V))\to C^\infty(Gr_i(V))$$
by
\begin{eqnarray}\label{E:Int1}
(T_\alp f)(E)=\int_{F\in Gr_i(V)}|\cos(E,F)|^\alp f(F) dF,
\end{eqnarray}
where $dF$ is the $O(n)$-invariant Haar probability measure on $Gr_i(V)$.
\begin{remark}\label{Int-R:2}
One can show that the integral (\ref{E:Int1}) absolutely converges for $Re(\alp)>-1$ (see  \cite[Lemma 2.1]{alesker-alpha}).
\end{remark}
It is well known (see below) that $T_\alp$ has a meromorphic continuation in $\alp\in \CC$.
For $\alp_0\in \CC$ we will denote by $S_{\alp_0}$ the
first non-zero coefficient in the decomposition of the meromorphic function $T_\alp$ near $\alp_0$, namely
$$T_\alp=(\alp-\alp_0)^k\cdot( S_{\alp_0} +O(\alp-\alp_0))\mbox{ as } \alp\to \alp_0, \, S_{\alp_0}\ne 0.$$
Thus $$S_{\alp_0}\colon C^\infty(Gr_i(V))\to C^\infty(Gr_i(V))$$
is defined and non-zero for any $\alp_0\in \CC$, and coincides with $T_{\alp_0}$ for all but countably many $\alp_0$
(necessarily $Re(\alp_0)\leq -1$ for exceptional values of $\alp_0$). $S_\alp$ will also be called the $\alp$-cosine transform.

\hfill

The $\alp$-cosine transform on Grassmannians, including the case of even functions on the sphere (e.g. $i=1$),
was studied by both analysts and geometers for a long period of time: \cite{koldobsky,ournycheva-rubin, rubin-inversion,rubin-adv,rubin-funk-cosine, schneider,semjanistyi}.
The case $\alp=1$ plays a special role in convex and stochastic geometry, see \cite{goodey-howard1,goodey-howard2,goodey-howard-reeder}, and in particular in valuation theory
\cite{alesker-bernstein,alesker-jdg-03,alesker-gafa-lefschetz,alesker-fourier}.

\hfill

Obviously the $\alp$-cosine transform $S_\alp$ (and also $T_\alp$) commutes with the natural action of the orthogonal group $O(n)$ on functions on
the Grassmannian. In representation theory there is a detailed information available about the action of $O(n)$ and $SO(n)$ on functions on Grassmannians
(see Section \ref{S:diff-oper} below). This allows one to apply tools from the representation theory of $O(n)$ and $SO(n)$ to the problem. This method was used  in
\cite{goodey-howard1,goodey-howard2,goodey-howard-reeder} in the case $\alp=1$.

\hfill

It was observed in \cite{alesker-bernstein} for $\alp=1$ and in \cite{alesker-alpha} for general $\alp$ that the $\alp$-cosine transform
can be rewritten in such a way to commute with an action of the much larger full linear group $GL(V)$. In that language it is a map
$$S_\alp\colon C^\infty(Gr_i(V),L_\alp)\to C^\infty(Gr_{n-i}(V),M_\alp),$$
where $L_\alp,M_\alp$ are $GL(V)$-equivariant line bundles over Grassmannians, and a choice of a Euclidean metric on $V$ induces
identification $Gr_{n-i}(V)\simeq Gr_i(V)$ by taking the orthogonal complement, and $O(n)$-equivariant identifications of $L_\alp$ and $M_\alp$
with trivial line bundles. For the explicit description of $L_\alp, M_\alp$ see (\ref{E:P01})-(\ref{E:P02}) in Section \ref{S:invariant-cosine} below.
Moreover it was observed in \cite{alesker-alpha}
that the $\alp$-cosine transform is essentially a special case of the well known representation theoretical construction of intertwining integrals
(see e.g. \cite{wallach-book}). These observations opened a way to use the infinite dimensional representation theory of $GL(V)$ in the study of the $\alp$-cosine transform.
Indeed in some more recent studies of the $\alp$-cosine transform \cite{genkai-zhang,olafsson-pasquale,gourevitch-sahi2} the representation theory of $GL(V)$
played the key role.

\hfill

Very recently the second and the third authors \cite{gourevitch-sahi-unique} have proven the following
new representation theoretical characterization of the $\alp$-cosine and Radon transforms. Let $L$ and $M$ be $GL(V)$-equivariant
complex line bundles over real Grassmannians $Gr_i(V)$ and $Gr_j(V)$ respectively. Then the space of $GL(V)$-equivariant
linear continuous maps $$C^\infty(Gr_i(V),L)\to C^\infty(Gr_j(V),M)$$ is at most 1-dimensional. Moreover
if $L$ and $M$ are $O(n)$-equivariantly isomorphic to trivial bundles, then the only such maps are either
Radon transform\footnote{The non-trivial observation that the Radon transform between Grassmannians can be rewritten in $GL(V)$-equivariant terms was first made in \cite{gelfand-graev-rosu}; see also Section \ref{S:invariant-radon} of this paper for the details.}, or the $\alp$-cosine transform $S_\alp$ for some $\alp\in \CC$, or a rank one operator. The last case
is degenerate and less interesting. Nevertheless it can be described explicitly as follows. All such operators
are obtained by a twist by a character of $GL(V)$ from a single operator mapping smooth measures on $Gr_i(V)$ to smooth functions on $Gr_j(V)$ given by first integrating a measure and then imbedding this number to functions as a constant function. These three cases are not quite independent: for example in Theorem \ref{T:cosine-radon}  of this paper we will show that in the special case when $j=n-i$ the Radon transform from $Gr_i(V)$ to $Gr_{n-i}(V)$ is proportional to the $\alp$-cosine transform for $\alp=-\min\{i,n-i\}$.

\hfill

The goal of the present paper is to obtain an analytic description of $S_\alp$. We use representation theoretical tools of both $GL(n,\RR)$ and $O(n)$.
Let us describe our main results more precisely. To formulate our first main result let us fix again a Euclidean metric on $V$. Denote $r:=\min\{i,n-i\}$.
Then $S_\alp\colon C^\infty(Gr_r(V))\to C^\infty(Gr_r(V))$. In Section \ref{S:diff-oper}, Theorem \ref{T:D-nu}, we write down  explicitly
an $O(n)$-equivariant differential operator on functions on $Gr_r(V)$, denoted by $\hat\cd_{\nu}$ for any $\nu\in \CC$, which satisfies the following
property.
\begin{theorem}\label{T:main1}
If $\alp\not\in [-(r+1),-2]\cap \ZZ$ then for some constant $c_\alp\in \CC$ one has
$$S_\alp=c_\alp\cdot \hat \cd_{\frac{\alp}{2}}\circ S_{\alp+2}.$$
\end{theorem}
The formula for $\hat\cd_{\frac{\alp}{2}}$ is explicit but somewhat too
technical to be presented in the introduction.
\begin{example}\label{Example-sphere}
For $r=1$ the operator $\hat\cd_{\frac{\alp}{2}}$ is a linear combination of  the Laplace-Beltrami operator on the unit sphere
(restricted to even functions on the sphere) and the identity operator (see Example \ref{Ex:sphere-case} in Section \ref{S:diff-oper} below). For $\alp=-1$ on the sphere this operator was written down explicitly
in \cite[Proposition 2.1]{goodey-weil}.
Notice also that for $r=1$ the operator $S_{-1}$ is proportional to the Radon transform (this fact is well known and seems to be a folklore).
\end{example}
\begin{remark}
We have the following mutually excluding cases.

1) For $Re(\alp)>-1$ $S_\alp=T_\alp$ is given by the explicit formula (\ref{E:Int1}).

2) If $Re(\alp)\leq -1$, but $\alp\not\in -\NN$, then by repeated use of Theorem \ref{T:main1} the operator $S_\alp$
can be expressed as a composition of some explicit (though complicated) $O(n)$-equivariant differential operator and $S_\beta$
with $Re(\beta) >-1$ (thus $S_\beta$ is given by (\ref{E:Int1})). Clearly $Re(\beta)$ can be chosen to be arbitrarily large.

3) If $\alp\in -r-2\NN$ then $S_\alp$ is composition of an explicit $O(n)$-equivariant differential operator with $S_{-r}$.

4) If $\alp\in -(r+1)-2\NN$ then $S_\alp$ is a composition of an explicit $O(n)$-equivariant differential operator with $S_{-(r+1)}$.

5) It remains to consider $\alp=-(r+1),-r,\dots,-1$. (Notice however that by Theorem \ref{T:main1} one has
$S_{-1}=const\cdot \hat \cd_{-\frac{1}{2}}\circ S_1$.) Our second main result below is an analytic interpretation of
$S_\alp$ in all these cases except $\alp=-(r+1)$ for which we do not know a good interpretation.
\end{remark}

\begin{theorem}\label{T:main2}
Let $r:=\min\{i,n-i\}$ as previously. We have

1) $S_{-r}\colon C^\infty(Gr_i(V))\to C^\infty(Gr_{n-i}(V))$ is the Radon transform; it can be rewritten as a $GL(V)$-equivariant
operator (see Section \ref{S:invariant-radon} for details; if $i=n/2$ the Radon transform coincides with the identity operator by convention).

2) Let $\alp=-(r-1),\dots,-2,-1$. Then  $S_{\alp}\colon C^\infty(Gr_i(V))\to C^\infty(Gr_{n-i}(V))$
can be presented in two different ways:

(a) as a composition of two Radon transforms via the Grassmannian $Gr_{|\alp|}(V)$;

(b) as a composition of two Radon transforms via the Grassmannian $Gr_{n-|\alp|}(V)$.

In both cases (a) and (b) one can rewrite the spaces and the operators in $GL(V)$-eqivariant language (see Section \ref{S:composition-Radon} for details).
\end{theorem}
This result is a combination of Theorems \ref{T:cosine-radon} and \ref{P:compos} below.
\begin{example}
Let $r=1$. This is the case of even functions on the sphere. As we have already mentioned in Example \ref{Example-sphere} above,
Theorem \ref{T:main2} says that $S_{-1}$ is (up to a constant) the Radon transform on the sphere.
\end{example}

\hfill

%*******************************************************************

The last main result of this paper is Theorem \ref{T:main3} below which computes explicitly
one more characteristic of the $\alp$-cosine transform $S_\alp$: the support of the distributional kernel.
We refer for the details to Section \ref{S:support}. Here we remind only that by the Schwartz kernel theorem
$S_\alp$ is given by a kernel on $Gr_i(V)\times Gr_{n-i}(V)$ which is a generalized section of an appropriate line bundle
over that space.  Since $S_\alp$ is $GL(V)$-equivariant operator, this section is $GL(V)$-invariant. Hence its support is a compact $GL(V)$-invariant
subset of $Gr_i(V)\times Gr_{n-i}(V)$.

To describe this support, let us observe that $GL(V)$ has finitely many orbits on $Gr_i(V)\times Gr_{n-i}(V)$. Each orbit consists of
pairs of subspaces $(E^i,F^{n-i})$ such that the dimension of their intersection is equal to a given number. Thus for each $l=0,1,\dots,r$ define the orbit
$$\co_l:=\{(E^i,F^{n-i})|\, \dim (E^i\cap F^{n-i})=l\}.$$
Clearly $\bar\co_l\supset \co_k$ iff $l\leq k$.

\begin{theorem}\label{T:main3}
(a) If $\alp\ne -1,-2,-3,-4\dots$ then the support of the
distributional kernel is $\bar \co_0$, i.e. maximal.

(b) If $\alp=-1,-2,\dots,-r+1$, then the support is equal
$\bar{\co}_{|\alp|}$.

(c) If $\alp \in -r-2\ZZ_{\geq 0}$ then the support is $\bar\co_r$, i.e. minimal.

(d) If $\alp\in -r-1-2\ZZ_{\geq 0}$ then the support is $\bar
\co_{r-1}$ (i.e. next to minimal).
\end{theorem}

\hfill

{\bf Acknowledgements.} The first author thanks F. Gonzalez and B. Rubin for several useful correspondences.
We thank Thomas Wannerer for a careful reading of the first version of the paper and numerous remarks.

\section{Some linear algebraic isomorphisms.}\label{S:linear-alg}
In this section we describe several canonical isomorphisms from (multi-) linear algebra.
All this material is well known and standard, but we would like to review it for the sake of completeness. In the rest of the paper
we will use these isomorphisms even
without mentioning them explicitly.

Let $X$ be an $n$-dimensional real vector space. In this section all constructed spaces will be real, but in subsequent
sections we will sometimes use the same notation for real spaces (more precisely, real vector bundles) and their complexifications. We denote
$$\det X:=\wedge^n X.$$
Let $Y\subset X$ be an $m$-dimensional linear subspace. Then we have a canonical isomorphism
\begin{eqnarray}\label{E:Iso1}
\det Y\otimes \det (X/Y)\tilde\to \det X
\end{eqnarray}
which is given by
$$(v_1\wedge\dots\wedge v_l)\otimes(\xi_1\wedge\dots\wedge \xi_{n-l})\mapsto v_1\wedge\dots\wedge v_l\wedge\bar\xi_1\wedge\dots\wedge\bar\xi_{n-l},$$
where $\bar \xi_i\in V$ is a lift of $\xi_i\in X/Y$. It is easy to see that the map (\ref{E:Iso1}) is well defined and is an isomorphism.

We have another canonical isomorphism
\begin{eqnarray}\label{E:Iso2}
\det X^*\simeq (\det X)^*
\end{eqnarray}
which is induced by the perfect pairing $\det X^*\times \det X\to\RR$ given by
$$(v_1^*\wedge\dots\wedge v_n^*,w_1\wedge\dots\wedge w_n)\mapsto \sum_{\sigma}sgn(\sigma)\prod_{i=1}^nv_i^*(w_{\sigma(i)}),$$
where the sum runs over all the permutations of length $n$.

\hfill

Let us define the orientation space $or(X)$ of the vector space $X$. Let $\cb(X)$ denote the set of all basis of $X$.
The group $GL(X)$ acts freely and transitively on $\cb(X)$. By definition $or(X)$ is the space of functions $f\colon \cb(X)\to \RR$
such that
$$f(g(e))=sgn(\det(g))f(e)$$
for any basis $e\in \cb(X)$ and any $g\in GL(X)$. Clearly $or(X)$ is 1-dimensional.

Any basis of $X$ defines a basis of $or(X)$ which is the function from $\cb(X)$ equal to 1 on this particular basis.
Let us define a canonical isomorphism
\begin{eqnarray}\label{E:Iso3}
or(Y)\otimes or(X/Y)\tilde\to or(X).
\end{eqnarray}
Choose a basis $v_1,\dots,v_m\in Y$; let $f_Y\in or(Y)$ be the function from $\cb(Y)$ equal to 1 on this basis. Fix a basis $\xi_1,\dots,\xi_{n-m}\in X/Y$; let
$f_{X/Y}\in or(X/Y)$ be the function equal to 1 on this basis.  Fix a lift $\bar \xi_i\in X$
for $i=1,\dots,n-m$. Then $v_1,\dots v_m,\bar\xi_1,\dots,\bar\xi_{n-m}$ is a basis of $X$. Let $f_X$ be the function equal to 1 on that basis.
Then there exists unique isomorphism (\ref{E:Iso3}) such that
$$f_Y\otimes f_{X/Y}\mapsto f_X.$$ It is not hard to check that this map is independent of the choices of bases.

By taking dual basis one easily constructs canonical isomorphisms
\begin{eqnarray}\label{E:Iso4}
or(X^*)\simeq or(X)\simeq or(X)^*.
\end{eqnarray}

Let $D(X)$ denote the space of real valued Lebesgue measures on $X$ (as we mentioned at the beginning of this section, in subsequent sections
$D(X)$ will denote the complexification of this space, namely the space of complex valued Lebesgue measures.) $D(X)$ is 1-dimensional.
There is a canonical isomorphism
\begin{eqnarray}\label{E:Iso5}
\det X^*\otimes or(X)\tilde\to D(X)
\end{eqnarray}
which is described as follows. Fix a basis $e_1,\dots,e_n\in X$. Let $f\in \cb(X)$ be the function equal to 1 on this basis.
Let $e_1^*,\dots,e_n^*\in X^*$ be the dual basis. Then there exists a unique isomorphism (\ref{E:Iso6}) such that
$(e_1^*\wedge\dots\wedge e_n^*)\otimes f$ is mapped to the Lebesgue measure whose value on the parallelepiped spanned by $e_1,\dots,e_n$
is equal to 1.

Isomorphisms (\ref{E:Iso1}), (\ref{E:Iso2}), (\ref{E:Iso3}), (\ref{E:Iso5}) imply the canonical isomorphisms
\begin{eqnarray}\label{E:Iso6}
D(Y)\otimes D(X/Y)\simeq D(X),\\
D(X^*)\simeq (D(X))^*.
\end{eqnarray}
One can easily check the following canonical isomorphisms
\begin{eqnarray}\label{E:Iso6.55}
D(X\otimes Y)\simeq D(X)^{\otimes \dim Y}\otimes D(Y)^{\otimes\dim X},\\
\det(X\otimes Y)\simeq \det(X)^{\otimes \dim Y}\otimes \det(Y)^{\otimes\dim X},\\
or(X\otimes Y)\simeq or(X)^{\otimes \dim Y}\otimes or(Y)^{\otimes\dim X}
\end{eqnarray}

\hfill

For $\alp\in \CC$ let us define the space $D^\alp(X)$ of $\alp$-densities on a real $n$-dimensional vector space $X$ as
the space of functions $f\colon \cb(X)\to \RR$ satisfying
$$f(g(e))=|\det g|^\alp f(e) \mbox{ for any } g\in GL(X),\, e\in \cb(X).$$
Clearly $\dim D^\alp(X)=1$.
We have a canonical isomorphism
$$D(X)\tilde\to D^1(X).$$

It is defined as follows: any Lebesgue measure $\mu$ on $X$ defines the function on $\cb(X)$ whose value
at a basis $e\in \cb(X)$ is equal to the $\mu$-measure of the parallelepiped spanned by the basis $e$.

Furthermore we have the following canonical isomorphisms (where $Y\subset X$ is a linear subspace):
\begin{eqnarray*}
D^\alp(Y)\otimes D^\alp(X/Y)\simeq D^\alp(X),\\
D^\alp(X)\otimes D^\beta(X)\simeq D^{\alp+\beta}(X),\\
(D^\alp(X))^*\simeq D^\alp(X^*)\simeq D^{-\alp}(X).
\end{eqnarray*}

\section{Invariant form of the $\alpha$-cosine transform.}\label{S:invariant-cosine}
Let us describe an invariant description of the $\alpha$-cosine transform following \cite{alesker-alpha}.
This is going to be a $GL(V)$-invariant operator
$$T_\alp\colon C^\infty(Gr_i(V),L_\alp)\to C^\infty(Gr_{n-i}(V),M_\alp),$$
where $L_\alp,M_\alp$ are $GL(V)$-equivariant line bundles over corresponding Grassmannians, $n=\dim V$.
The line bundles are described as follows. For subspaces $E^i\in Gr_i(V),\, F^{n-i}\in Gr_{n-i}(V)$ the fibers are
\begin{eqnarray*}
L_\alp|_{E^i}=D^\alp(V/E^i)\otimes|\omega_{Gr_i}|\big|_{E^i},\quad
M_\alp|_{F^{n-i}}=D^\alp(F^{n-i}),
\end{eqnarray*}
where $D^\alp(X)$ denotes the space of $\alpha$-densities on a vector space $X$, and $|\omega_{Gr_i}|$ denotes the line bundle of densities over $Gr_i(V)$, i.e. its fiber over a point $p$
is equal to the space of complex valued Lebesgue measures on the tangent space of $p$. In this notation we have
\begin{eqnarray}\label{E:001}
(T_\alp f)(F)=\int_{E\in Gr_i(V)} pr^*_{F\to V/E}(f(E)),
\end{eqnarray}
where $pr_{F\to V/E}\colon F\to V/E$ is the restriction to $F$ of the canonical projection $V\to V/E$, and the star $*$ denotes the pull-back on $\alpha$-densities.
It is obvious that $T_\alp$ is well defined for $Re(\alp)\geq 0$, namely the integral (\ref{E:001}) converges. Actually by \cite{alesker-alpha}, Lemma 1.1, the integral
converges absolutely for $Re(\alp)>-1$. Since $T_\alp$ coincides with the intertwining integral, by a general result (see e.g. \cite[Theorem 10.1.6]{wallach-book}) $T_\alp$ has a meromorphic continuation in $\alp\in \CC$.
For any $\alp_0\in \CC$ we will denote by $S_{\alp_0}$ the operator which is the leading non-zero coefficient in the decomposition of $T_\alp$ at $\alp_0$. Thus
\begin{eqnarray}\label{E:P001}
S_{\alp}\colon C^\infty(Gr_i(V),L_\alp)\to C^\infty(Gr_{n-i}(V),M_\alp)
\end{eqnarray}
is a non-zero $GL(V)$-equivariant operator.

For future reference, let us rewrite the line bundles $L_\alp, M_\alp$ slightly more explicitly.
%We will use the following notation.
%Let $l$ be a real oriented line, and let $\gamma\in \CC$ be a complex number. Let $l_{>0}$ denote the positive open half-line.
%The dual line $l^*$ admits a unique orientation such that $(l^*_{>0},l_{>0})>0$. Then denote by $l^{\otimes\gamma}$ the space of $\CC$-valued
%functions on $l^*_{>0}$ which are $\gamma$-homogeneous. Clearly  $l^{\otimes\gamma}$ is a complex line. If $\gamma=k$ is an integer, this line is canonically isomorphic
%to the complexification of $l^{\otimes k}$ is the usual sense. Moreover it is easy to see that for a vector space $W$
%$$D^{-\gamma}(W)=|\det (W)|^{\otimes \gamma}.$$
Below for a vector space $W$ we denote $D^{-\gamma}(W)$ by $|\det(W)|^{\otimes\gamma}$. Then we have
\begin{align}\label{E:P01}
&L_\alp|_{E^i}=|\det E^i|^{\otimes (n+\alp)}\otimes |\det V|^{\otimes -(i+\alp)},\\\label{E:P02}
&M_\alp|_{F^{n-i}}=|\det F^{n-i}|^{\otimes -\alp},
\end{align}
where we have use the well known canonical isomorphism for the tangent space
$$T_{E^i}Gr_i(V)=(E^i)^*\otimes V/E^i.$$

\section{Invariant form of the Radon transform.}\label{S:invariant-radon}

Let us denote by $R_{ji}$ the Radon transform from the $i$-Grassmannian to the $j$-Grassmannian in an $n$-dimensional vector space $V$. In order
to describe it in $GL(V)$-invariant terms it is convenient to consider two cases: $i<j$ and $i>j$, while for $i=j$ the Radon transform $R_{ii}$ is the identity operator.
First the Radon transform between Grassmannians was rewritten in $GL(V)$-equivariant  term in \cite{gelfand-graev-rosu}.

\hfill

\underline{\bf{Case $i<j$.}} In this case
$$R_{ji}\colon C^\infty(Gr_i(V),L)\to C^\infty(Gr_j(V),M),$$
where
\begin{align*}
&L|_{E^i}=|\det E^i|^{\otimes j} \mbox{ for any } E^i\in Gr_i(V),\\
&M|_{F^j}=|\det F^j|^{\otimes i} \mbox{ for any } F^j\in Gr_j(V).
\end{align*}
Then $$(R_{ji}f)(F^j)=\int_{E^i\in Gr_i(F^j)}f(E^i),$$
where the last expression makes sense because
$$f(E^i)\in L|_{E^i}= M|_{F^j}\otimes |\omega_{Gr_{i}(F^j)}|\big|_{E^i},$$
where $|\omega_{Gr_{i}(F^j)}|$ denotes the line bundle of densities on $Gr_i(F^j)$ so that its sections
can be integrated over $Gr_{i}(F^j)$.

\hfill

\underline{\bf{Case $i>j$.}} In this case
$$R_{ji}\colon C^\infty(Gr_i(V),S)\to C^\infty(Gr_j(V),T),$$
where
\begin{align*}
&S|_{E^i}=|\det E^i|^{\otimes (n-j)}\otimes |\det V|^{\otimes -(i-j)},\\
&T|_{F^j}=|\det F^j|^{\otimes (n-i)}.
\end{align*}
Moreover for $E^i\supset F^j$
$$S|_{E^i}=T|_{F^j}\otimes |\omega|\big|_{E^i},$$
where $|\omega|$ denotes the line bundle of densities over the manifold of
$i$-dimensional linear subspaces containing the given $j$-dimensional subspace $F^j$.
Hence the invariant form of $R_{ji}$ is given by
$$(R_{ji}f)(F^j)=\int_{E^i\supset F^j}f(E^i).$$

\hfill

It was shown in \cite{gelfand-graev-rosu} that $R_{ji}$ has maximal possible range, more precisely $R_{ji}$
is onto if $\dim Gr_i(V)\geq \dim Gr_j(V)$, and it is a closed imbedding if $\dim Gr_i(V)\leq \dim Gr_j(V)$. In particular
$R_{ji}$ is an isomorphism if $j=n-i$.

\begin{theorem}\label{T:cosine-radon}
(1) Let $i\ne n/2$, $r:=\min\{i,n-i\}$. Then the $(-r)$-cosine transform $S_{-r}$ is proportional to the Radon transform $R_{n-i,i}$.

(2) Let $i=n/2$, thus $r=n/2$. Then the cosine transform $S_{-r}$ is proportional to the identity operator on the space of sections of
the line bundle over $Gr_{\frac{n}{2}}(V)$ whose fiber over $E\in Gr_{\frac{n}{2}}(V)$ is equal to $|\det E|^{\otimes \frac{n}{2}}$.
\end{theorem}
{\bf Proof.} In the case (1) $S_{-r}$ and $R_{n-i,i}$ act between the same pair of line bundles by the above discussion of this section and
by (\ref{E:P01})-(\ref{E:P02}). Hence by the uniqueness theorem \cite{gourevitch-sahi-unique},
which says that the space of $GL(V)$-equivariant
operators between spaces of sections of $GL(V)$-equivariant line bundles over two Grassmannians is at most one dimensional\footnote{In this paper the uniqueness theorem \cite{gourevitch-sahi-unique} is used in a special
case of a pair of $i$- and $(n-i)$-dimensional Grassmannians. In this case it easily follows from Theorem
\ref{T:12} of this paper.}, they must be proportional.
In the case (2) $S_{-r}$ acts on the space of sections of the line bundle described in the statement of the theorem due to (\ref{E:P01})-(\ref{E:P02}).
Hence by the uniqueness theorem of \cite{gourevitch-sahi-unique} again $S_{-r}$ must be proportional to the identity operator.

\section{Composition of Radon transforms}\label{S:composition-Radon}
In either case of the previous section ($i<j$ or $i>j$), if we twist the source and the target spaces by the same character of $GL(V)$
and twist the map $R_{ji}$ accordingly, then the obtained map also will be $GL(V)$-equivariant, and it is natural
to call it also Radon transform. We will denote by $R'_{ji}$ such a twist of the Radon transform $R_{ji}$ without specifying
the character. Now we will discuss when two (twisted) Radon transforms are composable as $GL(V)$-equivariant operators, namely
$R'_{kj}\circ R'_{ji}$ is well defined. It is easy to see that necessarily $k=n-i$ and there are two cases: $i<j>n-i$ or $i>j<n-i$.
Moreover, using a global twist by a character of $GL(V)$,
we may assume that one of the composed twisted Radon transforms is not twisted in fact.

\underline{Case 1: $i<j>n-i$.} We have a composition of maps
\begin{eqnarray*}
C^\infty(Gr_i(V),L)\overset{R_{ji}}{\to}C^\infty(Gr_j(V),M)\overset{R_{n-i,j}'}{\to}C^\infty(Gr_{n-i},T),
\end{eqnarray*}
where
\begin{align}\label{E:R1}
&L|_{E^i}=|\det E^i|^{\otimes j} \mbox{ for any } E^i\in Gr_i(V),\\\label{E:R2}
&M|_{F^j}=|\det F^j|^{\otimes i} \mbox{ for any } F^j\in Gr_j(V),\\\label{E:R3}
&T|_{G^{n-i}}=|\det G^{n-i}|^{\otimes (n-j)}\otimes |\det V|^{\otimes (i+j-n)} \mbox{ for any } G^{n-i}\in Gr_{n-i}(V).
\end{align}

\hfill

\underline{Case 2: $i>j<n-i$.}  We have a composition of maps
\begin{eqnarray*}
C^\infty(Gr_i(V),L)\overset{R_{ji}}{\to}C^\infty(Gr_j(V),M)\overset{R_{n-i,j}'}{\to}C^\infty(Gr_{n-i},T),
\end{eqnarray*}
where
\begin{align*}
&L|_{E^i}=|\det E^i|^{\otimes (n-j)},\\
&M|_{F^j}=|\det F^j|^{\otimes (n-i)}\otimes |\det V|^{\otimes (i-j)},\\
&T|_{G^{n-i}}=|\det G^{n-i}|^{\otimes j}\otimes |\det V|^{\otimes (i-j)}.
\end{align*}

\hfill

Now let us observe that $i>j<n-i$ if and only of $i<n-j>n-i$. Hence according to cases 1 and 2, for $j<\min\{i,n-i\}$
we have two compositions of the Radon transforms between the same spaces:
\begin{eqnarray}\label{E:p1}
R'_{n-i,j}\circ R_{ji}, R'_{n-i,n-j}\circ R_{n-j,i}\colon C^\infty(Gr_i(V),L)\to C^\infty(Gr_{n-i}(V),T),\\\label{E:p2}
\mbox{ where } L|_{E^i}=|\det E^i|^{\otimes (n-j)},\, T|_{G^{n-i}}=|\det G^{n-i}|^{\otimes j}\otimes |\det V|^{\otimes (i-j)}.
\end{eqnarray}

Now let us observe that, due to (\ref{E:P01})-(\ref{E:P02}) and (\ref{E:p1}-\ref{E:p2}), the $(-j)$-cosine transform twisted by a character acts between the same spaces:
$$S_{-j}\otimes |\det (\cdot)|^{\otimes(i-j)}\colon  C^\infty(Gr_i(V),L)\to C^\infty(Gr_{n-i}(V),T).$$

\begin{theorem}\label{P:compos}
Let $1\leq j<\min\{i,n-i\}$. Then the three $GL(V)$-intertwining operators $R'_{n-i,j}\circ R_{ji}, R'_{n-i,n-j}\circ R_{n-j,i},$ and $S_{-j}\otimes |\det (\cdot)|^{\otimes(i-j)}$
between $$C^\infty(Gr_i(V),L)\to C^\infty(Gr_{n-i}(V),T)$$
are proportional to each other.
\end{theorem}
{\bf Proof.} This immediately follows from the uniqueness theorem \cite{gourevitch-sahi-unique}. \qed

\section{The precise statement and proof of Theorem \ref{T:main1}.}\label{S:diff-oper}
In this section we will use representation theory of $SO(n)$ and $O(n)$ on functions on real Grassmannians.
First let us remind few standard representation theoretical facts (see e.g. \cite{zhelobenko}).

All irreducible (continuous) representations of $SO(n)$ or $O(n)$ are necessarily finite dimensional.
For $SO(n)$ they can be parameterized by highest weights. In turn, highest weights of $SO(n)$ can be identified
with some combinatorial data, namely with sequences of integers
$$(m_1,\dots,m_{l-1},m_l) \mbox{ where } l:=\lfloor n/2\rfloor,$$
which satisfy
\begin{eqnarray*}
m_1\geq\dots m_{l-1}\geq m_l\geq 0 \mbox{ if } n \mbox{ is odd},\\
m_1\geq\dots m_{l-1}\geq |m_l|  \mbox{ if } n \mbox{ is even}.
\end{eqnarray*}
The natural representation of $SO(n)$ or $O(n)$ on $L^2(Gr_i(\RR^n))$ is isomorphic to such a representation
on $L^2(Gr_{n-i}(\RR^n))$ using the $O(n)$-equivariant identification $Gr_i(\RR^n)\simeq Gr_{n-i}(\RR^n)$ by taking the orthogonal complement.
Hence denote $r:=\min\{i,n-i\}$ and consider the action of $SO(n)$ and $O(n)$ on $L^2(Gr_r(\RR^n))$. These are unitary representations and each irreducible
representation of either group is known to enter $L^2(Gr_r(\RR^n))$ with multiplicity at most 1. Moreover each $SO(n)$- or $O(n)$-irreducible subspace
consists of infinitely smooth functions.

It is well known that an irreducible $SO(n)$-module with the highest weights $(m_1,\dots,m_l)$ does enter the decomposition of $L^2(Gr_r(\RR^n))$
if and only if:

$\bullet$ all $m_i$ are even integers;

$\bullet$ $m_i=0$ for $i>r$.

\hfill

In addition we will need a description of the decomposition of $L^2(Gr_r(\RR^n))$ under the action of $O(n)$.
Let us denote for $1\leq d\leq r\leq n/2$
\begin{align*}
&\Lambda_r:=\{(m_1,\dots,m_r)\in (2\ZZ)^r|\, m_1\geq \dots m_r\geq 0\},\\
&\Lambda_{d,r}:=\{ m\in \Lambda_r|\, m_d=\dots =m_r=0.\}.
\end{align*}
Let us consider
two cases: $r<n/2$ and $r=n/2$.

\underline{Case 1.} $r<n/2$.

Each $SO(n)$-irreducible subspace of $L^2(Gr_r(\RR^n))$ is $O(n)$-invariant (and of course irreducible). Hence
we can write
$$L^2(Gr_r(\RR^n))=\oplus_{ m\in \Lambda_r} \ch_{ m},$$
where $\ch_{ m}$ is an irreducible representation of $O(n)$ such that its restriction to $SO(n)$
is still irreducible and has highest weight $(m_1,\dots,m_r,\underset{l-r \mbox{ times}}{\underbrace{0,\dots,0}})$.

\hfill

\underline{Case 2.} $r=n/2$.

Each $SO(n)$-irreducible subspace of $L^2(Gr_r(\RR^n))$ with highest weight $\bar m:=(m_1,\dots,m_{r-1},m_r=0)$ is $O(n)$-invariant
(and clearly irreducible). We will denote this subspace $\ch_{m}$.

However if $m_r\ne 0$ the $SO(n)$-irreducible subspace is {\itshape not} $O(n)$-invariant.
What happens it is that the sum of two $SO(n)$-irreducible subspaces with highest weights $(m_1,\dots,m_{r-1},m_r)$ and
$(m_1,\dots,m_{r-1},-m_r)$ is $O(n)$-invariant and $O(n)$-irreducible subspace. Let us denote this subspace by $\ch_{ m}$
where $ m:=(m_1,\dots,m_{r-1},|m_r|)$. Thus for $r=n/2$ we have again
$$L^2(Gr_r(\RR^n))=\oplus_{ m\in \Lambda_r}\ch_{ m}.$$

\hfill

To summarize, in both cases ($r<n/2$ and $r=n/2$) we obtain that
$$L^2(Gr_r(\RR^n))=\oplus_{ m\in \Lambda_r}\ch_{ m},$$
where $\ch_{ m}$ are pairwise non-isomorphic irreducible representations of $O(n)$ described above.

This language will be very useful
in our paper due to Theorem \ref{genkai-zhang-thm} below. To formulate it,
let us observe that any (suitably continuous) $O(n)$-equivariant operator on functions on a Grassmannians $Gr_r(\RR^n)$, in particular
the $\alp$-cosine transforms, leaves each $\ch_{ m}$ invariant. By the Schur's lemma it acts on $\ch_{ m}$ by multiplication
by a scalar.

Let us denote by $c_{\nu,r}(m)$ this scalar for the
$2\nu$-cosine transform on Grassmannian
manifold of rank $r$ (i.e. on $Gr_r(V)$ or $Gr_{n-r}(V)$). We denote by
$(\nu)_k:=\nu(\nu+1)\dots(\nu+k-1)$ for $k\in \ZZ_{\geq 0}$ and $\nu\in \CC$
(where for $k=0$ one defines $(\nu)_0:=1$).
\begin{theorem}[\cite{genkai-zhang}, see also \cite{olafsson-pasquale}]\label{genkai-zhang-thm}
\begin{eqnarray}\label{E:ap-1}
c_{\nu,r}( m)=N_\nu
\prod_{j=1}^r\frac{(\nu+\frac{j+1}{2}-\frac{m_j}{2})_{\frac{m_j}{2}}}
{(\nu+\frac{n}{2}-\frac{j-1}{2})_{\frac{m_j}{2}}},
\end{eqnarray} where $N_\nu$
is a meromorphic function of $\nu$ which can be written explicitly.
\end{theorem}

\begin{remark}\label{R:0.5}
Let us observe that zeros of numerators in the expression for
$c_{\nu,r}(m)$ are disjoint from zeros of the denominators.
Indeed, zeros of the numerator are of the form $-(j-1)/2 +x$, where
$j=1,\dots,r$, and $x\in \ZZ_{\geq 0}$, while zeros of the
denominator are of the form $-n/2+(k-1)/2-y$, where again
$k=1,\dots,r$, and $y\in \ZZ_{\geq 0}$. The smallest zero of the
numerator is $-(r-1)/2$, and the largest zero of denominator is
$-n/2+(r-1)/2$. Clearly $-n/2+(r-1)/2 <-(r-1)/2$.
\end{remark}

Let us denote
$$c'_{\nu,r}( m):=\frac{c_{\nu,r}( m)}{N_\nu}.$$
Let $T'_{2\nu}:=\frac{1}{N_\nu}T_{2\nu}$; it has the eigenvalues
$c'_{\nu,r}( m)$.
\begin{lemma}\label{L:1}
(1) $T'_{2\nu}$ depends meromorphically on $\nu\in \CC$. It has no
zeros.

(2) The poles of $T'_{2\nu}$ are precisely at $\nu\in
-\frac{n-j+1}{2}-\ZZ_{\geq 0},$ where $j=1,\dots,r$.

(3) The multiplicity $\mu(l)$ of the pole of $T'_{2\nu}$ at $\nu=l\in
\frac{1}{2}\ZZ$ is computed as follows:

(a) if $l>-\frac{n}{2}+\frac{r-1}{2}$ then $\mu(l)=0$;

(b) if $l\leq -\frac{n}{2}$ then: if $l+n/2\in \ZZ$ then
$\mu(l)=\lceil\frac{r}{2}\rceil$; and if $l+n/2\not\in \ZZ$ then
$\mu(l)=\lfloor\frac{r}{2}\rfloor$;

(c) if $-\frac{n}{2}<l\leq -\frac{n}{2}+\frac{r-1}{2}$  then
$\mu(l)=\lfloor\frac{r-1-n}{2}-l\rfloor+1$.
\end{lemma}
{\bf Proof.} (1), (2) follow from Remark \ref{R:0.5}; the absence of
zeros follows from the observation that $c'_{\nu,r}(0)=1$.
Part (3) follows by counting zeros of denominators in the right hand side of (\ref{E:ap-1}). \qed

\hfill

Let us define
\begin{eqnarray}\label{E:Diff}
d_{\nu,r}( m)=\frac{c'_{\nu,r}( m)}{c'_{\nu+1,r}(
m)}=\prod_{j=1}^r\frac{(\nu+\frac{j+1}{2}-\frac{m_j}{2})\cdot
(\nu+\frac{n}{2}-\frac{j-1}{2}+\frac{m_j}{2})}{(\nu+\frac{j+1}{2})\cdot(\nu+\frac{n}{2}-\frac{j-1}{2})}.
\end{eqnarray}

Let us denote by $\cd_\nu$ the operator with eigenvalues
$d_{\nu,r}( m)$. We will see below that $\cd_\nu$ is a differential
operator with coefficients depending meromorphically on $\nu$.
Clearly
\begin{eqnarray}\label{E:recursion}
T'_{2\nu}=\cd_\nu\circ T'_{2\nu+2}.
\end{eqnarray}

Clearly (\ref{E:Diff}) immediately implies the following lemma.
\begin{lemma}\label{L:D-operator}
The operator $\cd_\nu$ does not vanish identically for any $\nu$,
and the poles are precisely at $\nu=-\frac{j+1}{2} \mbox{ and } -\frac{n}{2}+\frac{j-1}{2} \mbox{ with } j=1,\dots,r$
(notice that if $r=n/2$ then $-\frac{r+1}{2}$ appears twice in the list). Moreover if $r\ne n/2$ then all the poles are simple;
if $r=n/2$ then the pole
at $\nu=-\frac{r+1}{2}$ has multiplicity 2, while all the other poles are simple.
\end{lemma}

%Let us choose an entire non-vanishing function $\phi(\nu)$ which has the same poles, counting multiplicities,
%as $T'_\nu$ (this is possible say due to the Mittag-Leffler theorem, see \cite{shabat-book}, \S 13, Section 42).

Clearly for any $\nu_0\in \CC$ the operator $S_{2\nu_0}$  is equal to
$$ S_{2\nu_0}=\lim_{\nu\to \nu_0}(\nu-\nu_0)^\mu\cdot T'_{2\nu},$$
where $\mu$ is the multiplicity of the pole of $T'_{2\nu}$ at $\nu_0$.

Similarly we denote the %differential
operator
$$\tilde \cd_{\nu_0}:=\lim_{\nu\to \nu_0}(\nu-\nu_0)^\tau\cd_{\nu},$$
where $\tau$ is the multiplicity of the pole of $\cd_{\nu}$ at $\nu_0$.

Let us also denote by $\hat\cd_{\nu}$ the operator
$$\hat\cd_{\nu}:=\left(\prod_{j=1}^r(\nu+\frac{j+1}{2})\cdot(\nu+\frac{n}{2}-\frac{j-1}{2})\right)\cdot \cd_{\nu}.$$
From (\ref{E:Diff}) it is clear that for any $\nu\in \CC$ the operator $\hat\cd_\nu$ acts on the irreducible $O(n)$-representation $\ch_m$ as multiplication by a scalar
$$\prod_{j=1}^r(\nu+\frac{j+1}{2}-\frac{m_j}{2})\cdot
(\nu+\frac{n}{2}-\frac{j-1}{2}+\frac{m_j}{2}).$$
Thus it is clear that for any $\nu$ the operator $\hat\cd_\nu$ does nos not vanish identically, and $\hat \cd_\nu$ and $\tilde\cd_\nu$ are proportional
to each other. Below we will see, following \cite{itoh}, that $\hat\cd_\nu$ is an $O(n)$-invariant differential operator on functions on the Grassmannian which will be written down
more explicitly.

Our next goal is the following result.
\begin{proposition}\label{P:exceptions}
Let $l\not\in [-\frac{r+1}{2},\dots,-1]\cap\frac{1}{2}\ZZ$.
Then there exists a constant $c_l\in \CC$ such that
$$S_{2l}=c_l\cdot\hat\cd_l\circ S_{2l+2}.$$
\end{proposition}
{\bf Proof.} The proof is by a careful investigation of the orders of the poles in (\ref{E:recursion}). Below we will denote by
$\mu(T'_{2l})$ (resp. $\mu(\cd_l)$) the order of the pole of $T'_{2l}$ (resp. $\cd_l$) at $l$.

Since the poles $\nu$ of $T'_{2\nu}$ and $\cd_\nu$ are half-integers, the result obviously holds for $l$ not a half-integer.
Thus below we will always assume that {\itshape $l\in \frac{1}{2}\ZZ$} and $l$ does not belong to the segment $[-\frac{r+1}{2},\dots,-1]$.

\hfill

\underline{Case 1: $l>-1$.}  In this case $\mu(T'_{2l})=\mu(T'_{2l+2})=\mu(\cd_l)=0$. Hence in this case
the result follows immediately from (\ref{E:recursion}).

\hfill

\underline{Case 2: $l\leq -\frac{n}{2}-1$.} In this case $\mu(T'_{2l})=\mu(T'_{2l+2})=:\mu$ and both are equal either to
$\lceil r/2\rceil$ or $\lfloor r/2 \rfloor$ by Lemma \ref{L:1}(3)(b). Thus $\mu(\cd_l)=0$. Multiplying both sides
of (\ref{E:recursion}) by $(\nu-l)^\mu$ we conclude the proposition is this case.

\hfill

\underline{Case 3: $l=-\frac{n}{2}-\frac{1}{2}$.} We have by Lemma \ref{L:1}
$$ \mu(T'_{2l})=\mu(T'_{2l+2})=\lfloor r/2\rfloor, \mu(\cd_l)=0.$$
Multiplying (\ref{E:recursion}) by $(\nu-l)^{\lfloor r/2\rfloor}$ we deduce the result.

\hfill

\underline{Case 4: $l=-\frac{n}{2}$.} We have
\begin{eqnarray}\label{E:case 4}
\mu(\cd_l)=1, \mu(T'_{2l})=\lceil r/2 \rceil.
\end{eqnarray}

\underline{Subcase 4a: assume $r\geq 3$.} Then $l+1\in (-\frac{n}{2},-\frac{n}{2}+\frac{r-1}{2}]$. Hence
$$\mu(T'_{2l+2})=\lfloor\frac{r-1}{2}\rfloor.$$
It is easy to see that $\lceil r/2 \rceil=\lfloor\frac{r-1}{2}\rfloor +1$. Hence $ \mu(T'_{2l})=\mu(T'_{2l+2})+1$.
Hence multiplying (\ref{E:recursion}) by $(\nu-l)^{\mu(T'_{2l})}$ be get the proposition in subcase 4a.

\underline{Subcase 4b: assume $r=1,2$.} By Lemma \ref{L:1}(3)(b)
$$\mu(T'_{2l+2})=0.$$
But by (\ref{E:case 4})
$$\mu(T'_{2l})=\mu(\cd_l)=1.$$
Multiplying (\ref{E:recursion}) by $\nu-l$ we get the result in the subcase 4b. Thus case 4 is proven.

\hfill

\underline{Case 5: $-\frac{n}{2}<l\leq -\frac{n}{2}+\frac{r-1}{2}$.} First recall that from the very beginning we have
assumed that $l\not\in [-\frac{r+1}{2},\dots,-1]$. Hence if $l=-\frac{n}{2}+\frac{r-1}{2}$, in that case it follows that $r<n/2$.
Lemma  \ref{L:1}(3) implies that
$$\mu(T'_{2l})=\mu(T'_{2l+2})+1.$$
Also $\mu(\cd_l)=1$ (here we have used $r<n/2$ if $l=-\frac{n}{2}+\frac{r-1}{2}$; otherwise the order of the pole would be equal to 2).
Hence multiplying (\ref{E:recursion}) by $(\nu-l)^{\mu(T'_{2l})}$ we get the proposition
in case 5.

\hfill

\underline{Case 6: $-\frac{n}{2}+\frac{r-1}{2}<l<-\frac{r+1}{2}$.} This is the last case to be checked.
In this case we necessarily have $r< n/2$. We have
$$\mu(\cd_l)=0.$$
Also by Lemma \ref{L:1}(3)(a)
$$\mu(T'_{2l})=\mu(T'_{2l+2})=0.$$
Hence the result follows. \qed

\hfill

%***************************************************************

Now we are going to describe $\hat\cd_\nu$ as a differential operator. Recall that the Pfaffian of a
$2k\times2k$ skew-symmetric matrix $M=\left(  M_{ij}\right)  $ is
\begin{equation}
Pf\left(  M\right)  =\frac{1}{k!}\sum_{\sigma}sgn\left(  \sigma\right)  M_{\sigma\left(
1\right)  \sigma\left(  2\right)  }\cdots M_{\sigma\left(  2k-1\right)
\sigma\left(  2k\right)  } \label{=Pf}%
\end{equation}
where the sum is over all permutations $\sigma$ satisfying $\sigma\left(
1\right)  <\sigma\left(  2\right)  ,\ldots,\sigma\left(  2k-1\right)
<\sigma\left(  2k\right)  $. We note that this definition makes sense even if
the entries of $M$ belong to a non-commutative algebra, such as an enveloping algebra.

Now let $X$ be the $n\times n$ skew-symmetric matrix whose $ij$-th entry is%
\[
X_{ij}=E_{ij}-E_{ji}\in\mathfrak{o}\left(  n\right)  \text{,}
\]
where $E_{ij}$ is the $n\times n$ matrix all whose entries vanish but the entry on $i$th raw and $j$th column
is equal to 1.

For $I\subset\left\{  1,\ldots,n\right\}  $ let $X_{I}$ be the principal
$I$-minor of $X$; if $\left\vert I\right\vert $ is even then the Pfaffian
$Pf\left(  X_{I}\right)  ,$ defined by (\ref{=Pf}), is an element of the
enveloping algebra $U=U\left(  \mathfrak{o}\left(  n\right)  \right)  $. Observe that in the sum  (\ref{=Pf})
for $Pf(X_I)$ in each summand all terms commute with each other, while terms from different summands may not commute.

For
$d\leq n/2$ we define
\begin{eqnarray}\label{E:V-d}
V_{d}=\left(  -1\right)  ^{d}\sum_{\left\vert I\right\vert =2d}Pf\left(
X_{I}\right)  ^{2}\in U(\mathfrak{o}(n)).
\end{eqnarray}
\begin{theorem}[\cite{itoh-umeda}]
$V_{d}$ belongs to $\left(U\left(  \mathfrak{o}\left(  n\right)
\right)\right)^{O(n)}  $.
%and vanishes on $K$-types with highest weights in $\Lambda_{d,r}$
\end{theorem}

\begin{example}
(1) $V_0=1\in U(\mathfrak{o}(n))$.

(2) $V_1$ is proportional to the Casimir element of $\mathfrak{o}(n)$.
\end{example}

\begin{remark}
It was shown in \cite[Theorem 2.3]{gonzalez-kakehi} that the algebra
$(U(\mathfrak{o}(n)))^{O(n)}$ is a polynomial algebra on (independent) generators $V_1,V_2,\dots,V_{\lfloor n/2\rfloor}$.
However we will not use this fact in the paper.
\end{remark}

\begin{theorem}\label{T:D-nu}
We have $\hat\cd_\nu=(-\frac{1}{4})^r\sum_{k=0}^rc_kV_k$ where
$$c_k=\prod_{j=k+1}^r[j+2\nu +1][j-(2\nu+n+1)].$$
\end{theorem}
Notice that in the last theorem $V_k\in (U(\mathfrak{o}(n)))^{O(n)}$ are identified with the differential operators
they induce on the Grassmannian $Gr_r(\RR^n)$.

\begin{example}\label{Ex:sphere-case}
Clearly $V_0$ induces the identity operator on functions on all Grassmannians. Since $V_1\in U(\mathfrak{o}(n))$ is proportional to the Casimir element,
it induces on functions on each Grassmannian an operator proportional to the Laplace-Beltrami operator. This if $r=1$ then $\hat \cd_\nu$ is a linear combination
of the Laplace-Beltrami and the identity operators.
\end{example}

To prove Theorem \ref{T:D-nu} we will need some preparations.
\begin{lemma}\label{L:vanishing}
For $1\leq d\leq r\leq n/2$ the action of the element $V_d$ vanishes on irreducible $O(n)$-modules $\ch_m$ with
$m\in \Lambda_{d,r}$.
\end{lemma}
{\bf Proof.} The irreducible representations of $O(n)$ from $\Lambda_{d,r}$ are precisely those appearing in
$L^2(Gr_{d-1}(\RR^n))$; each has a vector invariant under $\mathfrak{m}:=\mathfrak{o}(d-1)\times\mathfrak{o}(n-d+1)$ (the so called spherical vector).
Let us show that each summand in the definition of $Pf(X_I)$, $|I|=2d$, (see (\ref{=Pf})) has at least one
factor from $\mathfrak{m}$; this obviously will imply the lemma.

Let us assume the contrary. Then there exists a subset $$I=\{i_1<i_2<\dots<i_{2d}\}\subset\{1,\dots,n\}$$ and a permutation $\sigma$ of length $2d$ such that
$X_{i_{\sigma(1)}i_{\sigma(2)}},\dots, X_{i_{\sigma(2d-1)}i_{\sigma(2d)}}\not\in \mathfrak{m}.$ This implies that for the odd indices
$$1\leq i_{\sigma(1)}, i_{\sigma(3)},\dots,i_{\sigma(2d-1)}\leq d-1.$$
Since all these $d$ numbers are distinct, we get a contradiction. \qed

\hfill

Let $1\leq r\leq n/2$. Let $\mathbb{D=D}\left(  Gr_r(\RR^n)\right)  $ be
the algebra of $O\left(  n\right)  $-invariant differential operators on
$ Gr_r(\RR^n) $, and let $\mathbb{D}_{k}$ be the subspace of operators of order $\leq
k$. Also let $\mathbb{P=P}\left(  z_{1},\ldots,z_{r}\right)  $ be the algebra
of polynomials in $\left(  z_{1},\ldots,z_{r}\right)  $, which are invariant
under sign changes and permutations of coordinates, namely under all transformations of the form
$$(z_1,\dots,z_r)\mapsto (\pm z_{\sigma(1)},\dots,\pm z_{\sigma(r)})$$ with $\sigma$ being any permutation of length $r$.
Let $\mathbb{P}_{k}$ be the subspace
of polynomials with degree $\leq k$. Harish-Chandra's theorem for the
symmetric space $ Gr_r(\RR^n)$ can be formulated as follows.

\begin{theorem}
\label{HC}There is a unique algebra isomorphism $\gamma:\mathbb{D\tilde\rightarrow P}$,
such that it maps $\mathbb{D}_{k}$ onto $\mathbb{P}_{k}$, and $D$ $\in
\mathbb{D}$ acts on the irreducible $O(n)$-module $\ch_\mu $ by the scalar
$\gamma\left(  D\right)  \left(  \tilde{\mu}\right)  $, where%
\[
\tilde{\mu}:=\mu+\rho,\rho:=\left(  \rho_{1},\ldots,\rho_{r}\right)  ,\rho
_{j}:=n/2-j.
\]

\end{theorem}
We failed to find a proof of this theorem in literature; its proof will be given in the appendix below.

\begin{lemma}\label{L:sahi6}
Let $1\leq d\leq r\leq n/2$. Let $\rho=(\rho_1,\dots,\rho_r)\in \CC^r$ be an arbitrary vector. The space
$$\cv_d:=\{p\in\PP_{2d}|\, p(\mu+\rho)=0 \mbox{ for all } \mu\in \Lambda_{d,r}\}$$
is at most one dimensional.
\end{lemma}
{\bf Proof.} Fix $d\geq 1$ and prove the lemma by the induction in $r\geq d$.

\underline{Step 1.} Assume $r=d$.

\underline{Case 1.} Assume $\rho_d=0$. Then $p(\mu_1,\dots,\mu_{d-1},0)=0$ for any $\mu_1,\dots,\mu_{d-1}\in \CC$.
Hence $p$ is divisible by $z_d$. Since $p$ is invariant under sign changes, it is divisible by $z_d^2$. Since $p$ is invariant
under permutations, it is divisible by $z_1^2\dots z_d^2$. But $\deg p\leq 2d$. Hence $p$ is proportional to $z_1^2\dots z_d^2$,
and the lemma follows in this case.

\underline{Case 2.} Assume $\rho_d\ne 0$. Then $p(\mu_1,\dots,\mu_{d-1},\rho_d)=0$ for any $\mu_1,\dots,\mu_{d-1}\in \CC$.
Hence $p$ is divisible by $z_d-\rho_d$ for any $\mu_1,\dots,\mu_{d-1}\in \CC$. But since $p$ is invariant under sign changes and all
permutations, it is divisible by $\prod_{i=1}^r(z_i^2-\rho_i^2)$. Since $\deg p\leq 2d$, $p$ is proportional to the latter polynomial.

\hfill

\underline{Step 2.} Assume now $r>d$.
Let $q$ be another polynomial in $z_1,\dots,z_r$ satisfying the same assumptions of the lemma as $p$. We have to show that $p$ and $q$ are
proportional.

Define the polynomial
$$\hat p(z_1,\dots,z_{r-1}):=p(z_1,\dots,z_{r-1},\rho_r).$$
Clearly $\hat p$ is invariant under all permutations and sign changes and $\deg \hat p\leq 2d$. Moreover $\hat p(\hat\mu+\hat \rho)=0$
for any $\hat\mu\in \Lambda_{d,r-1}$, where $\hat \rho=(\rho_1,\dots,\rho_{r-1})$. Similarly define $\hat q$ for $q$. By the induction assumption
$\hat p$ and $\hat q$ are proportional, say
$\hat p=\alpha \hat q.$ Define
$$\tau=p-\alp q.$$
We are going to show that $\tau \equiv 0$; this will imply the lemma. Clearly $\tau$ is a polynomial in $z_1,\dots,z_r$ of degree at most $2d$, invariant under permutations and sign
changes. Moreover
$$\hat \tau(z_1,\dots,z_{r-1}):=\tau(z_1,\dots,z_{r-1},\rho_r)$$
vanishes identically. As in Step 1 it follows that $\tau$ is divisible by $\prod_{i=1}^r(z_i^2-\rho_r^2)$ (no matter if $\rho_r$ vanishes or not).
But since $\deg \tau\leq 2d<2r$, it follows that $\tau\equiv 0$, i.e. $p=\alp q$. \qed

\hfill

Let us introduce more notation. Let $x:=(x_1,\dots,x_r)$. Let $e_{k}\left(  x\right)  $ and $h_{k}\left(  x\right)  $ be the elementary
and complete symmetric polynomials in $x$,
\[
e_{k}\left(  x\right)  :=\sum_{j_{1}<\cdots<j_{k}}x_{j_{1}}\cdots x_{j_{k}%
},\quad h_{k}\left(  x\right)  :=\sum_{j_{1}\leq\cdots\leq j_{k}}x_{j_{1}%
}\cdots x_{j_{k}},
\]
Then we have
\[
\prod_{i=1}^{r}\left(  1+tx_{i}\right)  =\sum_{k=0}^{r}e_{k}\left(  x\right)
t^{k}\text{\thinspace\quad}\prod_{i=1}^{r}\left(  1+tx_{i}\right)  ^{-1}%
=\sum_{k=0}^{\infty}\left(  -1\right)  ^{k}h_{k}\left(  x\right)  t^{k}%
\]

\begin{lemma}
\label{Sym}Let $e_{ij}\left(  x\right)  =e_{i}\left(  x_{j},\ldots
,x_{r}\right)  $ and $h_{ij}\left(  x\right)  =h_{i}\left(  x_{j},\ldots
,x_{r}\right)  $ then we have
\[
\sum_{j=i}^{k}\left(  -1\right)  ^{k-j}h_{k-j,k}e_{j-i,i+1}=\delta_{ik}.
\]

\end{lemma}

{\bf Proof.}
If $k=i$ then both sides are $1$, while if $k<i$ then both sides are $0$.
Finally, if $k>i$, then the left side is the coefficient $c_{k-i}$ of
$t^{k-i}$ in%
\[
\prod_{l=k}^{r}\left(  1+tx_{l}\right)  ^{-1}\prod_{l=i+1}^{r}\left(
1+tx_{l}\right)  =\prod_{l=i+1}^{k-1}\left(  1+tx_{l}\right)  .
\]
But this is a polynomial of degree $k-i-1$ in $t$, and hence $c_{k-i}=0$.
\qed

\begin{definition}
\label{def:ab}For $0\leq i,j\leq r$ we put $a_{ij}=b_{ij}=0$ if $i>j$, and
\[
a_{ij}=e_{j-i,i+1}\left(  x\right)  ,b_{ij}=\left(  -1\right)  ^{j-i}%
h_{j-i,j}\left(  x\right)  \text{ if }i\leq j.
\]

\end{definition}

\begin{corollary}
\label{ab}The matrices $\left(  a_{ij}\right)  $ and $\left(  b_{ij}\right)  $
are mutual inverses.
\end{corollary}

{\bf Proof.}
Lemma \ref{Sym} shows that $\sum_{j=0}^{r}a_{ij}b_{jk}=\delta_{ik}$, as desired.
\qed

Below we will always assume that $1\leq d\leq r\leq n/2$ and $$\rho_j=n/2-j,\, j=1,\dots,r.$$
Now the operator $V_d$ from (\ref{E:V-d}) and the isomorphism $\gamma$ from Theorem \ref{HC} satisfy
\begin{theorem}
We have $$\gamma\left(  V_{d}\right)  =\sum\nolimits_{k=0}^{d}\left(
-1\right)  ^{d-k}h_{d-k}\left(  \rho_{d}^{2},\ldots,\rho_{r}^{2}\right)
e_{k}\left(  z_{1}^{2},\ldots,z_{r}^{2}\right)  .$$
\end{theorem}

{\bf Proof.} Notice that for $d<n/2$ this is precisely Theorem 9.1(A) in \cite{kakehi-99}.
We will prove the theorem for $d=n/2$ (hence $r=n/2$) though this assumption will be used only in the last step in the computation of the constant of
proportionality.
We note that the expression on the right is the coefficient $v_{d}$ of $t^{d}$
in
\begin{equation}
\prod_{l=d}^{r}\left(  1+t\rho_{l}^{2}\right)  ^{-1}\prod_{l=1}^{r}\left(
1+tz_{l}^{2}\right)  . \label{=van}%
\end{equation}
We claim that $\gamma\left(  V_{d}\right)  $ and $v_{d}$ both belong to the
one-dimensional space $\mathcal{V}_{d}$. For $\gamma\left(  V_{d}\right)  $
this follows by Lemma \ref{L:vanishing}, while for $v_{d}$ we note that if
$z\in\left\{  \mu+\rho:\mu\in\Lambda_{d,r}\right\}  $ then
\[
z_{d}=\rho_{d},\ldots,z_{r}=\rho_{r}%
\]
Thus (\ref{=van}) is a polynomial of degree $<d$ in $t$ and hence
$v_{d}\left(  z\right)  =0$.

By Lemma \ref{L:sahi6} $\gamma\left(  V_{d}\right)  $ and $v_{d}$ are
proportional. It remains to show that the constant of proportionality is 1.
Now we are going to use the assumption $d=n/2$. It suffices to show that for some irreducible $O(n)$-module
$\ch_\mu$ one has $\gamma(V_d)(\mu+\rho)=v_d(\mu+\rho)\ne 0$. Under the assumption $d=n/2$ we have
$$V_d=(-1)^dPf(X)^2,$$
where $Pf(X)=\frac{1}{d!}\sum_{\sigma}'sgn(\sigma)X_{\sigma_1\sigma_2}\cdot\dots \cdot X_{\sigma_{n-1}\sigma_n}$ and the sum $\sum_{\sigma}'$
runs over all permutations $\sigma$ of length $n$ such that $\sigma_{2i-1}<\sigma_{2i}$ for all $i$. Since by assumption $d=r=n/2$, then $\rho_d=0$, and
the statement of the theorem reduces to equality
\begin{eqnarray}\label{E:equality-middle}
\gamma(V_d)=z_1^2\dots z_d^2.
\end{eqnarray}
We already know that the two polynomials are proportional, and we have to show that the constant of proportionality is equal to 1.
In order to prove that we will compute the action of $V_d$ on $\wedge^d \CC^{2d}$, where $\CC^{2d}$ is the standard representation of $O(2d)$.

$\wedge^d\CC^{2d}$ is an irreducible $O(2d)$-module (see e.g. \cite{fulton-harris}, Exercise 19.3). However as an $SO(2d)$-module it is a sum of two irreducible
$SO(2d)$-modules with highest weights $\mu_+=(\underset{d}{\underbrace{1,\dots,1,1}}),\, \mu_-=(\underset{d}{\underbrace{1,\dots,1,-1}})$ (see \cite{fulton-harris},
Theorem 19.2 and Remark (ii) on p. 289 there).
Adding $\rho$ to $\mu_{\pm}$ we get respectively in these cases
$$z_1=d,z_2=d-1,\dots,z_d=\pm 1.$$
Hence
\begin{eqnarray}\label{E:zz}
z_1^2 z_2^2\dots z_d^2=(d!)^2.
\end{eqnarray}
We are going to show that $V_d$ acts on $\wedge^d\CC^{2d}$ as $(d!)^2$; this will finish the proof of the theorem.
Let $e_1,\dots,e_{2d}$ be an orthonormal basis of $\CC^{2d}$. We have
\begin{align}
&Pf(X)e_1\wedge\dots\wedge e_d=\\\label{E:Pf-action}\frac{1}{d!}
&\sum_{\sigma} ' sgn(\sigma)\underset{X_{\sigma_1\sigma_2}}{\underbrace{(E_{\sigma_1\sigma_2}-E_{\sigma_2\sigma_1})}}
\dots\underset{X_{\sigma_{n-1}\sigma_n}}{\underbrace{(E_{\sigma_{n-1}\sigma_n}-E_{\sigma_n}E_{\sigma_{n-1}})}}e_1\wedge\dots\wedge e_d.
\end{align}
Clearly
\begin{eqnarray*}
E_{ij}e_p=\left\{\begin{array}{cc}
    e_i,&j=p\\
    0,&j\ne p
    \end{array} \right.
\end{eqnarray*}
Hence for $i\ne j$
\begin{eqnarray*}
X_{ij}e_p=\left\{\begin{array}{cc}
                   e_i,&j=p\\
                   -e_j,&i=p\\
                    0,&i,j\ne p
                 \end{array}\right.
\end{eqnarray*}
This implies that
$$X_{\sigma_{2i-1}\sigma_{2i}}X_{\sigma_{2l-1}\sigma_{2l}}e_p=0 \mbox{ for } i\ne l$$
since $\{\sigma_{2i-1},\sigma_{2i}\}\cap \{\sigma_{2l-1},\sigma_{2l}\}=\emptyset$.

Moreover if for some $1\leq p\leq d$ one has
$$X_{\sigma_{2i-1}\sigma_{2i}}e_p=\pm e_q \mbox{ with } 1\leq q\leq d,$$
then such a summand may be omitted from the sum (\ref{E:Pf-action}) since in this case
$\{\sigma_{2i-1},\sigma_{2i}\}=\{p,q\}$, and hence there is no another $X_{\sigma_{2l-1}\sigma_{2l}}$ to apply on $e_q$.

It follows that the only non-zero summands in (\ref{E:Pf-action}) correspond to $\sigma=(\sigma_1,\sigma_2,\dots,\sigma_{n-1},\sigma_n)$ of
the form:

$\bullet$ $\sigma_1,\sigma_3,\dots,\sigma_{n-1}$ is a permutation of $\{1,\dots,d\}$;

$\bullet$ $\sigma_2,\sigma_4,\dots,\sigma_{n}$ is a permutation of $\{d+1,\dots,n\}$.

Then
$$(\ref{E:Pf-action})=(-1)^d\sum_\tau sgn(1\tau_12\tau_2\dots d\tau_d)e_{\tau_1}\wedge\dots\wedge e_{\tau_d},$$
where the sum runs over all permutations $\tau=(\tau_1,\dots,\tau_d)$ of the set $\{d+1,\dots,n\}$.

Observe that
$$sgn(1\tau_12\tau_2\dots d\tau_d)=(-1)^{\frac{d(d-1)}{2}}sgn(\tau).$$
Hence
\begin{eqnarray}\label{E:Pf-action2}
Pf(X)e_1\wedge\dots\wedge e_d=(-1)^{d+\frac{d(d-1)}{2}}(d!)e_{d+1}\wedge\dots\wedge e_n.
\end{eqnarray}

\hfill

Now let us compute $Pf(X)e_{d+1}\wedge\dots\wedge e_n$. Since $\wedge^d \CC^{2d}$ is a multiplicity free $SO(2d)$-module and $Pf(X)$
and the Hodge star $\star$ commute with $SO(2d)$,\footnote{The fact that $Pf(X)$ commutes with $SO(2d)$ was proven in Proposition 3.6 in \cite{itoh-umeda}.}
it follows that $Pf(X)$ and $\star$ commute with each other. Then we have
\begin{multline}\label{E:Pf-action3}
Pf(X)e_{d+1}\wedge\dots\wedge e_n=Pf(X)\left(\star(e_1\wedge\dots\wedge e_d)\right)=
\star Pf(X)(e_1\wedge\dots\wedge e_d)\overset{(\ref{E:Pf-action2})}{=}\\
\star((-1)^{d+\frac{d(d-1)}{2}}(d!)e_{d+1}\wedge\dots\wedge e_n)=
(-1)^{d+\frac{d(d-1)}{2}}(d!)\cdot (-1)^de_{1}\wedge\dots\wedge e_d).
\end{multline}
Hence by (\ref{E:Pf-action2}) and (\ref{E:Pf-action3}) we get
\begin{eqnarray*}
V_d(e_1\wedge\dots\wedge e_d)=(-1)^dPf(X)^2e_1\wedge\dots\wedge e_d=\\
(d!)^2e_1\wedge\dots\wedge e_d\overset{(\ref{E:zz})}{=}\\
z_1^2\dots z_d^2\cdot e_1\wedge\dots\wedge e_d.
\end{eqnarray*}
The theorem is proved.
\qed

Let us write $\rho^{2}=\left(  \rho_{1}^{2},\ldots,\rho_{r}^{2}\right)  $,
then by Definition \ref{def:ab}, we have
\[
a_{ij}\left(  \rho^{2}\right)  =e_{j-i}\left(  \rho_{i+1}^{2},\ldots,\rho
_{r}^{2}\right)  ,b_{ij}\left(  \rho^{2}\right)  =\left(  -1\right)
^{j-i}h_{j-i}\left(  \rho_{j}^{2},\ldots,\rho_{r}^{2}\right)  ,\text{ for
}i\leq j.
\]
and so by the previous theorem, the eigenvalues of $V_{d}$ are
\[
\sum_{k=0}^{d}b_{kd}\left(  \rho^{2}\right)  e_{k}\left(  \tilde{\mu}_{1}%
^{2},\ldots,\tilde{\mu}_{r}^{2}\right)
\]
By Corollary \ref{ab}, we deduce the following result.

\begin{theorem}
\label{Ed}The operator $E_{d}=\sum_{k=0}^{d}a_{kd}\left(  \rho^{2}\right)
V_{k}$ has eigenvalues $e_{d}\left(  \tilde{\mu}_{1}^{2},\ldots,\tilde{\mu}%
_{r}^{2}\right)  $ on $\ch_\mu$.
\end{theorem}

{\bf Proof of Theorem \ref{T:D-nu}.} Recall that the eigenvalue of operator $\hat\cd_\nu$ on $\ch_\mu$ is equal
to
\begin{align*}
&  \prod_{j=1}^{r}\left[  \nu+\frac{j+1}{2}-\frac{\mu_{j}}{2}\right]  \left[
\nu+\frac{n}{2}-\frac{j+1}{2}+\frac{\mu_{j}}{2}\right]  \\
&  =\frac{1}{4^{r}}\prod_{j=1}^{r}\left[  \left(  2\nu+\frac{n}{2}+1\right)
^{2}-\tilde{\mu}_{j}^{2}\right]  =\left(  \frac{-1}{4}\right)  ^{r}\prod
_{j=1}^{r}\left(  \tilde{\mu}_{j}^{2}+\lambda\right)
\end{align*}
where $\lambda:=-\left(  2\nu+n/2+1\right)  ^{2}$. By Theorem \ref{Ed} we have
$$\hat\cd_\nu=\left(-\frac{1}{4}\right)^r\sum_{d=0}^r\lambda^{r-d}E_d=\left(-\frac{1}{4}\right)^r\sum_{k=0}^r c_k V_k,$$
where
\begin{eqnarray*}
c_{k}=\sum_{d=k}^{r}\lambda^{r-d}a_{kd}\left(  \rho^{2}\right)  =\sum
_{d=k}^{r}\lambda^{r-d}e_{d-k}\left(  \rho_{k+1}^{2},\ldots,\rho_{r}%
^{2}\right)  =\prod_{j=k+1}^{r}\left(  \rho_{j}^{2}+\lambda\right) =\\\prod_{j=k+1}^r[j+2\nu+1][j-(2\nu+n+1)] .
\end{eqnarray*}
\qed

\section{Support of the distributional kernel of the $\alp$-cosine transform.}\label{S:support}
The main goal of this section is to prove Theorem \ref{T:main3} of the introduction.
Let us fix an isomorphism $V\simeq \RR^n$. Let $e_1,\dots,e_n$ be the standard basis of $\RR^n$.
Let us consider open subset $\cu\subset Gr_{n-r}(\RR^n)$ consisting of subspaces $E^{n-r}$ trivially intersecting
the $r$-subspace $span\{e_{n-r+1},\dots,e_n\}$. Equivalently $E^{n-r}$ is equal to the span of columns of a matrix $\tilde X_{n\times (n-r)}$ of the form
\begin{eqnarray}\label{E:sup1}
\left[\begin{array}{c}
I_{n-r}\\
\hline
X
\end{array}\right],
\end{eqnarray}
where $X$ is an $r \times (n-r)$ matrix. Furthermore the map $E^{n-r}\mapsto X$ is a
diffeomorphism $\cu\tilde\to Mat_{r \times (n-r)}(\RR)$.

Let $\cl^\beta\to Gr_{n-r}(\RR^n)$ be the line bundle whose fiber over $E^{n-r}$ is equal to $\beta$-densities in $E^*$:
\begin{eqnarray*}
\cl^\beta|_{E^{n-r}}=|\det E|^{\beta}.
\end{eqnarray*}
Let us fix a trivialization of $\cl^\beta$ over $\cu$ as follows. Let $E^{n-r}\in \cu$ be an arbitrary subspace. Let $X\in Mat_{r\times(n- r)}(\RR)$ be the matrix corresponding to $E^{n-r}$.
Thus $E^{n-r}$ is the span of columns of the matrix $\tilde X$ defined by (\ref{E:sup1}). Let us denote the columns of $\tilde X$ by
$e_1(E),\dots,e_{n-r}(E)$. Then $e(E)^\beta:=|e_1(E)\wedge\dots\wedge e_{n-r}(E)|^{\otimes\beta}\in\cl^{\beta}|_{E^{n-r}}$ defines the required trivialization of
$\cl^\beta$ over $\cu$.

\hfill

Let $Q\subset GL_n(\RR)$ be the subgroup consisting of block diagonal matrices of the form $\left[\begin{array}{cc}
                                                                                                    A_{(n-r)\times(n-r)}&0\\
                                                                                                    0&B_{r\times r}
                                                                                                   \end{array}\right]$.

We are going to show that $\cu$ is $Q$-invariant under restriction of the natural action of $GL_n(\RR)$ on $Gr_{n-r}(\RR)$ to $\cu$, and compute the action of $Q$ on
$\cl^\beta$. We have for any $g=\left[\begin{array}{cc}
                                       A&0\\
                                       0&B
                                      \end{array}\right]\in Q$:
\begin{eqnarray*}
\left[\begin{array}{cc}
A&0\\
0&B
\end{array}\right]\left[\begin{array}{c}
                                I_{n-r}\\
                                X_{r \times (n-r)}
                                \end{array}\right]=\left[\begin{array}{c}
                                I_{n-r}\\
                                BXA^{-1}
                                \end{array}\right]\cdot A.
\end{eqnarray*}
From this equality we immediately conclude.
\begin{claim}\label{Cl:supAction}
1) $g$  maps subspace $E\in \cu$ corresponding to a matrix $X$ to subspace $E'\in\cu$ corresponding to
the matrix $BXA^{-1}$. In particular $\cu$ is $Q$-invariant.

2) $g$ maps $e(E)^\beta$ to $|\det A|^\beta\cdot e(E')^\beta$.
\end{claim}
Let us denote by $\cu_0$ the (closed) submanifold of $\cu$ consisting of subspaces containing $E_0^r:=span\{e_1,\dots,e_r\}$.
Clearly $\cu_0$ is identified with matrices of the form
$$\tilde X=\left[\begin{array}{c}
                  I_{n-r}\\
                  \hline
                  \begin{array}{c|c}
                  0_{r\times r}&Z
                  \end{array}
                  \end{array}\right],$$
where $Z$ is a $(n-2r)\times r$ matrix. Let us fix the transversal $\cn\subset \cu$ to $\cu_0$:
\begin{eqnarray}\label{E:sup2}
\cn=\left\{\left[\begin{array}{c}
                  I_{n-r}\\
                  \hline
                  \begin{array}{c|c}
                  Y_{r\times r}&0
                  \end{array}
                  \end{array}\right]\big| Y\in Mat_{r\times r}(\RR)\right\}\simeq Mat_{r\times r}(\RR),
\end{eqnarray}
where $0$ has size $(n-2r)\times r$.

Let us define the subgroup $Q_0\subset Q$ consisting of block-diagonal matrices of the form
$$\left[\begin{array}{ccc}
                  A_{r\times r}&0&0\\
                  0&I_{n-2r}&0\\
                  0&0&D_{r\times r}
                  \end{array}
                  \right].$$
The group $Q_0$ leaves $\cn$ invariant, and we have
\begin{eqnarray*}
\left[\begin{array}{ccc}
                  A_{r\times r}&0&0\\
                  0&I_{n-2r}&0\\
                  0&0&D_{r\times r}\end{array}
                  \right]\left[\begin{array}{c}
                  I_{n-r}\\
                  \hline
                  \begin{array}{c|c}
                  Y_{r\times r}&0
                  \end{array}
                  \end{array}\right]=
                  \left[\begin{array}{c}
                  I_{n-r}\\
                  \hline
                  \begin{array}{c|c}
                  DYA^{-1}&0
                  \end{array}
                  \end{array}\right]\left[\begin{array}{c|c}
                                             A&0\\
                                             \hline
                                             0&I_{n-2r}
                                             \end{array}
                                             \right].
\end{eqnarray*}
From the last equality we immediately conclude.
\begin{claim}\label{Cl:supAction2}
Under the identification $\cn\simeq Mat_{r\times r}(\RR)$ as in (\ref{E:sup2}),
an element $g=\left[\begin{array}{ccc}
                  A_{r\times r}&0&0\\
                  0&I_{n-2r}&0\\
                  0&0&D_{r\times r}\end{array}
                  \right]\in Q_0$ maps a subspace $E\in \cn$ corresponding to $Y$ to the subspace $E'\in \cn$
corresponding to $DYA^{-1}$, and maps $e(E)^\beta$ to $|\det A|^\beta\cdot e(E')^\beta$.
\end{claim}

\hfill

Now let us return back to the $\alp$-cosine transform. First consider the case $i\leq n-i$, thus $r=i$. Recall that we have the $\alp$-cosine transform
$$S_\alp\colon C^\infty(Gr_r(\RR^n), L_\alp)\to C^\infty(Gr_{n-r}(\RR^n),M_\alp),$$
where
\begin{eqnarray*}
L_\alp|_{E^r}=|\det E|^{\otimes(n+\alp)}\otimes|\det V|^{\otimes-(r+\alp)},\\
M_\alp|_{F^{n-r}}=|\det F^{n-r}|^{\otimes -\alp}.
\end{eqnarray*}

The distributional kernel of $S_\alp$ is a $GL_n(\RR)$-invariant generalized section $\mathbb{S}_\alp$ over $Gr_r(\RR^n)\times Gr_{n-r}(\RR^n)$
of the line bundle $(L_\alp^*\otimes|\omega_{Gr_r(\RR^n)}|)\boxtimes M_\alp$.
\begin{remark}
(1) We have also used the standard notation $\boxtimes$ of the exterior product of two bundles. Namely if $X_1,X_2$
are two manifold and $L_i,\, i=1,2,$ is a vector bundle over $X_i$ then by definition
$$L_1\boxtimes L_2:=p_1^*L_1\otimes p_2^*L_2,$$
where $p_i\colon X_1\times X_2\to X_i$ is the natural projection, $i=1,2$.

(2) Recall that $\mathbb{S}_\alp$ satisfies for any smooth section $f\in C^\infty(Gr_r(\RR^n),L_\alp)$
$$S_\alp(f)=\int_{Gr_r(\RR^n)\times Gr_{n-r}(\RR^n)}\mathbb{S}_\alp\cdot f,$$
and $\mathbb{S}_\alp$ is uniquely characterized by this property.
\end{remark}

Let $P$ denote the stabilizer of $E_0^r=span\{e_1,\dots,e_r\}$.
We will use the following theorem

\begin{theorem}[Frobenius descent, see e.g. {\cite[Theorem
4.2.3]{AGS1}}] \label{Frob}
Let a Lie group $K$ act on a smooth manifold $M$. Let $N$
be a smooth manifold with a transitive action of $K$. Let
$\phi:M \to N$ be a $K$-equivariant map.
Let $z \in N$ be a point and $M_z:= \phi^{-1}(z)$ be its fiber.
Let $K_z$ be the stabilizer of $z$ in $K$.
Let $\mathcal{E}$ be a $K$-equivariant vector bundle over $M$.

Then
there exists a canonical isomorphism between $K$-invariant generalized sections of $\mathcal{E}$ on $M$ and $K_z$-invariant generalized sections of $\mathcal{E}$ on $M_z$.
Moreover, for any closed $K$-invariant subset $Y\subset M$, the generalized sections supported in $Y$ are mapped to generalized sections supported in $Y\cap M_z$.
\end{theorem}

This isomorphism is the restriction to $M_z$. For the $K$-invariant generalized sections this restriction is well defined.

\begin{corollary}\label{L:sup3}
There is a natural one-to-one correspondence between $GL(V)$-invariant generalized sections of the line bundle $(L_\alp^*\otimes|\omega_{Gr_r(\RR^n)}|)\boxtimes M_\alp$
over $Gr_r(\RR^n)\times Gr_{n-r}(\RR^n)$ and $P$-invariant generalized sections of $M_\alp\otimes (L_\alp^*\otimes |\omega_{Gr_r(\RR^n)}|)|_{E^r_0}$ over $Gr_{n-r}(V)$.
Under this correspondence the support of the latter section is equal to the intersection of the support of the former one with the set $\{E_0^r\}\times Gr_{n-r}(V)$.
\end{corollary}

\begin{corollary}\label{L:sup4}
(1) There is a well defined restriction map from the space of $P$-invariant generalized sections of the line bundle
$M_\alp\otimes(L^*_\alp\otimes|\omega_{Gr_r(\RR^n)}|)\big|_{E_0^r}$ over $Gr_{n-r}(\RR^n)$ to the space of generalized functions
on $\cn\simeq Mat_{r\times r}(\RR)$. Generalized functions in the image of this map satisfy
\begin{eqnarray}\label{E:sup5}
S(AXB)=|\det A\cdot\det B|^{\alp}S(X)
\end{eqnarray}
for any  $X\in Mat_{r\times r}(\RR),\, A,B\in GL_r(\RR)$.

(2) This map is injective.

(3) The support of the restriction of a section is equal to the intersection of the support of that section with $\cn$.
\end{corollary}
{\bf Proof.} First restrict the generalized section to $\cu$. This restriction is injective since $\cu$ intersects all $P$-orbits in $Gr_{n-r}(\RR^n)$.
Now consider the map $\phi:\cu \to Mat_{r \times (n-2r)}(\RR)$ obtained by composing the diffeomorphism $\cu\tilde\to Mat_{r \times (n-r)}(\RR)$ with the operation of taking the right $n-2r$  columns. Let us define the subgroup $R\subset P$ consisting of matrices of the form
\begin{eqnarray}\label{E:def-g}
g=\left[\begin{array}{ccc}
                  A_{r\times r}&0&0\\
                  0&I_{n-2r}&0\\
                  0&E_{r\times n-2r}&D_{r\times r}
                  \end{array}
                  \right].
\end{eqnarray}
Note that $R$ preserves $\cu$. Let $R$ act on $Mat_{r \times (n-2r)}(\RR)$ by $g(Z)= DZ+E$, where $g\in R$ in given by (\ref{E:def-g}).
Note that this action is transitive, and the map $\phi:\cu \to Mat_{r \times (n-2r)}(\RR)$ is $R$-equivariant. The space $\cn$ is the fiber of the zero matrix under $\phi$, and the stabilizer in $R$ of $0\in Mat_{r \times (n-2r)}(\RR)$ is $Q_0$. The corollary follows now from the Frobenius descent.
\qed
\hfill

Now we are going to study generalized functions on $Mat_{r\times r}(\RR)$ satisfying (\ref{E:sup5}).
\begin{proposition}\label{P:1}
For any $\alp\in \CC$ the space of generalized functions $S$ satisfying (\ref{E:sup5})
is exactly one dimensional.
\end{proposition}

The positivity of the dimension follows from existence of
meromorphic continuation of the generalized function $|\det(\cdot)|^\alp$ (this is a special case of a general result due to J. Bernstein
\cite{bernstein-contin} which says that for any polynomial $P$ on $\RR^N$ the generalized function $|P|^\alp$ has a meromorphic continuation to the whole complex plane). To prove
uniqueness we will need some lemmas.
%\begin{lemma}\label{L:2}
%Let $0\leq k<r$. The tangent space to the submanifold of rank $k$
%matrices at the matrix $X_k:=\left[\begin{array}{cc}
%                         I_{k}&0\\
%                         0&0_{(r-k)\times (r-k)}
%                        \end{array}\right]$
%is equal to the space of matrices of the form
%$\left[\begin{array}{cc}
%                         A_{k\times k}&B_{k\times (r-k)}\\
%                         C_{(r-k)\times k}&0_{(r-k)\times (r-k)}
%                        \end{array}\right]$
%\end{lemma}
%with arbitrary matrices $A,B,C$.
%
%{\bf Proof.} The action of $GL_r\times GL_r$ on $Mat_{r\times r}(\RR)$
%is $(U,V)(X)=UXV^{-1}$. The corresponding action of the Lie algebra
%is $(U,V)(X)=UX-XV$. The tangent space to the rank $k$ matrices at
%$X_k$ is equal $UX_k-X_kV,\, U,V\in M_{r\times r}(\RR)$. A
%straightforward computation implies the result. \qed
%
%\hfill
%
%Lemma \ref{L:2} implies that
Define
$$\cn_k:=\left\{\left[\begin{array}{cc}
                      I_k&0\\
                      0&Y_{r-k,r-k}
                      \end{array}\right]\right\}.$$
%is a transversal to rank $k$ matrices. In fact it is easy to see
%that $\cn_k$ is transversal to all the strata it intersects (namely
%to all submanifolds of matrices of a fixed rank $\geq k$).

Note that using the Frobenius descent any generalized function $S$ satisfying (\ref{E:sup5}) on the open subset consisting of  matrices of rank at least $k$
can be
restricted to $\cn_k$. This restriction $\tilde S$ satisfies readily
the same equation (\ref{E:sup5}) but all matrices $A,B,X$ have size
$(r-k)\times (r-k)$.

\begin{lemma}\label{L:3}
If $S$ satisfies (\ref{E:sup5}) with parameter $\alp$ then its Fourier
transform satisfies (\ref{E:sup5}) with parameter $-r-\alp$.
\end{lemma}
The proof is straightforward.

\begin{lemma}\label{L:4}
If $S$ satisfies (\ref{E:sup5}) with parameter $\alp$ and is supported at 0 then $S$ is
proportional to
$$(\det(\pt_{ij}))^{2k}\delta(X),$$
where $\pt_{ij}$ is the partial derivative with respect to $x_{ij}$,
$k=0,1,2,\dots$, and $\alp=-r-2k\in -r-2\ZZ_{\geq 0}$.
\end{lemma}
{\bf Proof.} The Fourier transform $\FF(S)$ is a polynomial. By
Lemma \ref{L:3} it satisfies (\ref{E:sup5}) with parameter $-r-\alp$.
Hence $\FF(S)$ must be a polynomial of even degree. This implies the
lemma.\qed

\begin{proposition}\label{P:5}
Let $S\ne 0$ satisfy (\ref{E:sup5}) with parameter $\alp$. Let $0<l\leq r$
be an integer. Assume that $supp(S)=\{co-rank \geq l\}$. Then
$$\alp=-l, -l-2,-l-4, \dots.$$
\end{proposition}
{\bf Proof.} Let us restrict $S$ to $\cn_{r-l}$. This restriction
$\tilde S$ satisfies (\ref{E:sup5}) with the same $\alp$ but on the space
of matrices of size $l$. Hence our lemma follows from Lemma
\ref{L:4} applied to $\tilde S$. \qed

\hfill

We immediately deduce the following corollary.
\begin{corollary}\label{C:6}
If $S\ne 0$ satisfies (\ref{E:sup5}) with $\alp\ne -1,-2, -3,\dots$, then the
support of $S$ is equal to $Mat_{r\times r}(\RR)$. Moreover for such
values of $\alp$ the generalized function $S$ is unique up to a multiplicative constant.
\end{corollary}
Only uniqueness requires an explanation. On the open orbit the generalized
function is clearly unique up to a multiplicative constant.
Thus any two generalized functions $S$ might be
assumed to be equal on the open orbit. Then their difference is
supported on smaller orbits, and hence vanishes. This proves
Proposition \ref{P:1} in this case. Hence it remains to prove the
proposition only for $\alp=-1,-2,-3,\dots$. By applying the Fourier
transform and Lemma \ref{L:3} we see that it remains to prove the
proposition only for $\alp=-1,-2,\dots, -[\frac{r}{2}]$.
\hfill

\begin{lemma}\label{L:7}
The space of generalized functions $S$ satisfying (\ref{E:sup5}) with $\alp=-1$ is one
dimensional, and their support is equal to matrices of co-rank $\geq
1$ (exactly).
\end{lemma}
{\bf Proof.} Let us restrict any such generalized function to $\cn_{r-1}$. This
restriction is an even $(-1)$ -homogeneous generalized function on $\RR$. Hence
it is a multiple of the delta function $\delta(x)$. Let us show that
the kernel of this restriction map is zero. But the kernel consists
of generalized functions supported on matrices of co-rank 2 and higher. By
Proposition \ref{P:5} the parameter $\alp=-2,-3,-4,\dots$, i.e. not
$-1$; this is a contradiction. \qed

The above discussion shows that we proved Proposition \ref{P:1} for
$r=2,3$, namely we have:
\begin{corollary}\label{C:8}
Let $r=2,3$. Then for any $\alp$ the space of generalized functions satisfying
(\ref{E:sup5}) is one dimensional.
\end{corollary}

\hfill

{\bf Proof of Proposition \ref{P:1}.} We prove by induction by $r$.
For $r=2,3$ it is Corollary \ref{C:8}. Thus let us assume $r\geq 4$.
It remains to consider the case $\alp=-2,\dots, -[\frac{r}{2}]$.

Let us restrict $S$ to $\cn_1$. This restriction is a generalized function on
$Mat_{(r-1)\times (r-1)}(\RR)$ which satisfies (\ref{E:sup5}) with the
same $\alp$. By the induction assumption the space of such generalized functions is
one dimensional. Hence it remains to show that the kernel of the
restriction is zero. Indeed the kernel consists of generalized functions
supported at 0. Then by Lemma \ref{L:4}, $\alp\leq -r$. But
$-r<-[\frac{r}{2}]$. \qed

\begin{lemma}\label{L:9}
If $S\ne 0$ satisfies (\ref{E:sup5}) with $\alp=-r,-r-2,-r-4,\dots$ then
$$supp(S)=\{0\}.$$
\end{lemma}
{\bf Proof.} By the uniqueness part we see that $S$ is proportional
to $(\det (\pt_{ij}))^{2k}\delta(X)$, $ k=0,1,2\dots$. The result
follows. \qed

\begin{lemma}\label{L:10}
If $S\ne 0$ satisfies (\ref{E:sup5}) with $\alp=-1,-2,\dots,-r$. Then
$$supp (S)=\{co-rank\geq |\alp|\}.$$
\end{lemma}
{\bf Proof.} Let us restrict $S$ to $\cn_{r+\alp}\simeq
M_{|\alp|\times|\alp|}(\RR)$. $S$ must be proportional to $\delta(X)$. If
the coefficient of proportionality is not 0 then the lemma is proved.
Let us assume that it vanishes. Then let us choose $l$ such that
$supp(S)=\{co-rank\geq l\}$. Thus $l>|\alp|$. Then the restriction of
$S$ to $\cn_{r-l}\simeq M_{l\times l}(\RR)$ does not vanish and is
supported at 0. Hence $s=-l,-l-2,-l-4\dots$, which contradicts to the inequality $l>|\alp|$. \qed

\begin{lemma}\label{L:11}
If $S\ne 0$ satisfies (\ref{E:sup5}) with $\alp=-r-1,-r-3, -r-5,\dots$ then
$$supp(S)=\{rk\leq 1\}.$$
\end{lemma}
{\bf Proof.} Let us restrict $S$ to $\cn_1\simeq
M_{(r-1)\times(r-1)}(\RR)$. By Lemma \ref{L:9}
$supp(S|_{\cn_1})=\{0\}$. Equivalently $supp(S)\subset\{rk\leq 1\}$.
It remain to show that $supp(S)\ne \{0\}$. Indeed otherwise we can
apply Lemma \ref{L:4} and get a contradiction. \qed

\hfill

Let us summarize Proposition \ref{P:1}, Corollary \ref{C:6}, Lemma
\ref{L:9}, Lemma \ref{L:10}, Lemma \ref{L:11}. We get
\begin{theorem}\label{T:12}
(1) For any $\alp\in \CC$ the space of generalized functions $S$ satisfying
(\ref{E:sup5}) is exactly one dimensional.

(2) The support of $S$ is described completely by one of the
following cases:

(a) If $\alp\ne -1,-2,-3,\dots$ then $supp(S)=Mat_{r\times r}(\RR)$, i.e. is maximal.

(b) If $\alp=-1,-2,\dots,-r+1$, then $supp(S)=\{co-rank\geq |\alp|\}$.

(c) If $\alp\in -r-2\ZZ_{\geq 0}$, then $supp(S)=\{0\}$.

(d) If $\alp\in -r-1-2\ZZ_{\geq 0}$, then $supp(S)=\{rk\leq 1\}$.
\end{theorem}

Now Theorem \ref{T:main3} of the introduction follows immediately from
Lemma \ref{L:sup3}, Lemma \ref{L:sup4}, and Theorem \ref{T:12} in the case $r=i$, namely $i\leq n-i$.
It remains to consider the case $i>n-i$. This case easily follows from the previous one if one replaces
$E\in Gr_i(V)$ and $F\in Gr_{n-i}(V)$ by their annihilators $E^\perp\in Gr_{n-i}(V^*)$ and $F^\perp\in Gr_{i}(V^*)$.
Under this identification the $\alp$-cosine transform from $Gr_i(V)$ to $Gr_{n-i}(V)$ becomes the $\alp$-cosine transform
from $Gr_{n-i}(V^*)$ to $Gr_i(V^*)$. Hence  Theorem \ref{T:main3} is proved completely.

\appendix

\section{Appendix. Proof of Theorem \ref{HC}.}\label{S:appendix}
We start by recalling a few general basic facts on invariant differential operators on homogeneous spaces. Our main reference is \cite[Ch. II, \S 4.2]{gga},
where the case of connected groups is considered.
Let $G$ be a Lie group, let $\fg_0$ be the Lie algebra of $G$ and let $\fg:=\fg_0\otimes_\RR\CC$ be its complexification. Let $H\subset G$ be a closed subgroup
with Lie algebra $\fh_0$ and complexification $\fh$.

Let $D(G)$ denote the algebra of left $G$-invariant differential operators on $G$ which are invariant with respect to
all left translations. $D(G)$ is naturally isomorphic to the universal enveloping algebra
$U(\fg)$. Indeed the Lie algebra $\fg$ acts on functions on $G$ by
right infinitesimal translations; they commute with the left ones.
Hence we get a homomorphism of algebras $U(\fg)\to D(G)$. It is easy
to see that it is an isomorphism.

Let us denote by $D_H(G)$ the algebra of differential operators on $G$ which are left invariant under $G$ and right invariant under $H$.
Thus $D_H(G)\subset D(G)$. Under the above isomorphism $D(G)\simeq U(\fg)$, $D_H(G)$ is isomorphic to the subalgebra $U(\fg)^H$ of $H$-invariant elements.

Let $\pi\colon G\to G/H$ be the canonical map. The pull-back map is
$$\pi^*\colon C^\infty(G/H)\to C^\infty(G).$$
This gives a homomorphism of algebras
$$\tilde \pi\colon D_H(G)\to D(G/H)$$
which is uniquely characterized by the property
$\pi^*(\tilde\pi(D)(f))=D(\pi^*f)$ (here we have used that
$D(\pi^*f)$ is invariant under right translations by $H$, hence is a
pull-back under $\pi$ of a function from $C^\infty(G/H)$).

\begin{proposition}\label{P:homomor-operators}
Assume that there exists an $Ad(H)$-invariant subspace $\fm\subset
\fg$ such that $\fg=\fh\oplus\fm$ (for example this holds provided
$H$ is compact). Then the homomorphism $\tilde \pi\colon D_H(G)\to
D(G/H)$ is onto. Moreover its kernel is equal to $U(\fg)^H\cap
(U(\fg)\cdot \fh)$ under the above described isomorphism
$D_H(G)\simeq U(\fg)^H$.
\end{proposition}
{\bf Proof.} For a connected group $G$ the proposition was proved in
\cite{gga}, Ch. II, Theorem 4.6. We need to remove this assumption.
In fact we will not reduce the result to the connected case, but
give an independent proof.

First the algebra $\cd(X)$ of differential operators with smooth
coefficients on a manifold $X$ has a canonical filtration by the
differential operator's order $\{\cd(X)_{\leq k}\}$; it is compatible with the product. We
will prove first a more precise statement that $\tilde \pi$ is onto
on subspaces of operator of order at most $k$ for each $k$:
$$\tilde\pi\colon D_H(G)_{\leq k}\to D(G/H)_{\leq k}.$$
This is obviously true for $k=0$. Using induction it suffices to
show that $\tilde\pi$ is onto on all associated graded spaces
\begin{eqnarray}\label{E:assoc-graded}
\tilde\pi\colon D_H(G)_{\leq k}/D_H(G)_{\leq k-1}\to D(G/H)_{\leq
k}/D(G/H)_{\leq k-1}.
\end{eqnarray}
We have
\begin{eqnarray}\label{E:assoc1}
D_H(G)_{\leq k}/D_H(G)_{\leq k-1}\simeq \frac{U(\fg)^H_{\leq
k}}{U(\fg)^H_{\leq k-1}}=\left(\frac{U(\fg)_{\leq k}}{U(\fg)_{\leq
k-1}}\right)^H \simeq (Sym^k\fg)^H.
\end{eqnarray}
(Note that we have not used our assumption on $H$ in the above
isomorphisms.) On the other hand for any manifold $X$
$$\cd(X)_{\leq k}/\cd(X)_{\leq k-1}\simeq C^\infty(X,Sym^k (TX)),$$
where in the right hand side one has the space of smooth sections of
the $k$th symmetric power of the tangent bundle $TX$. Hence for
$X=G/H$ we get an imbedding
\begin{eqnarray}\label{E:assoc2}
D(G/H)_{\leq k}/D(G/H)_{\leq k-1}\inj (C^\infty(G/H,
Sym^k(T(G/H))))^G=(Sym^k (\fg/\fh))^H.
\end{eqnarray}
Thus using (\ref{E:assoc1}) and (\ref{E:assoc2}) we obtain from
(\ref{E:assoc-graded}) the maps
\begin{eqnarray}\label{E:assoc3}
(Sym^k\fg)^H\to D(G/H)_{\leq k}/D(G/H)_{\leq k-1}\inj(Sym^k
(\fg/\fh))^H.
\end{eqnarray}
Thus it suffices to prove that the composition of these two maps is
onto in order to show that the homomorphism $\tilde \pi$ is onto and
thus finish the proof of the first part of the proposition. But the
composed map is the canonical map
$$(Sym^k\fg)^H\to(Sym^k
(\fg/\fh))^H.$$ Its surjectivity follows immediately from the
assumption on existence of $Ad(H)$-invariant complement $\fm$ of
$\fh$ in $\fg$.

It remains to describe the kernel of $\tilde \pi\colon D_H(G)\to
D(G/H)$. By the Poincar\'e-Birkhoff-Witt theorem we have the
isomorphism of vector spaces
\begin{eqnarray}\label{E:PBW}
U(\fg)\simeq Sym^\bullet(\fm)\oplus U(\fg) \cdot \fh.
\end{eqnarray}
It is easy to see that an element $D\in U(\fg)$ considered as a left
$G$-invariant operator on functions on $G$ vanishes on all right
$H$-invariant functions if and only if under the isomorphism
(\ref{E:PBW}) it corresponds to an element of the form
$U(\fg) \cdot \fh$. \qed

\hfill

Now let us introduce some notation. For $0\leq p\leq q$ let us
denote
\begin{eqnarray*}
G_0=SO_0(p,q),\, K_0=SO(p)\times SO(q),\\
G_1=SO(p,q),\, K_1=S(O(p)\times O(q)),\\
G_2=O(p,q),\, K_2=O(p)\times O(q).
\end{eqnarray*}
Here $SO_0(p,q)$ denotes the connected component of the group
$SO(p,q)$.

\begin{proposition}\label{AP:1}
For any $i=0,1,2$ the imbedding $K_i\inj G_i$ induces an isomorphism on the groups of connected components.
\end{proposition}
{\bf Proof.} For $q=0$ the statement is trivial, and for $q=1$ it is
well known. Let us assume that $q\geq 2$ and proceed by induction on
$q$. For $i=0,1,2$ let us denote by $G_i'$ the group $G_i$ with
$(p,q-1)$ instead of $(p,q)$, and $K_i'$ the corresponding subgroup
of $G_i'$. It is easy to see that the standard imbedding $K_i'\inj
K_i$ induces an isomorphism on $\pi_0$. We have the commutative
diagram:
\begin{eqnarray}\label{Diag:connect}
\square[\pi_0(K_i')`\pi_0(G_i')`\pi_0(K_i)`\pi_0(G^i);```],
\end{eqnarray}
where the first horizontal map is an isomorphism by the induction assumption. It suffices to show that the second vertical arrow is an isomorphism.
Notice that for $i=0,1,2$
$$H:=G_i/G_i'=\{x_1^2+\dots+x_q^2-y_1^2-\dots-y_p^2=1\}.$$
Let us consider the smooth map $f\colon H\to \RR^p\times S^{q-1}$ given by
$$f(x_1,\dots,x_q,y_1,\dots,y_p)=\left((y_1,\dots,y_p),\frac{(x_1,\dots,x_q)}{(x_1^2+\dots +x_q^2)^{\frac{1}{2}}}\right).$$
Clearly $f$ is diffeomorphism. Hence for $q\geq 2$ the manifold $H$ is connected. From the exact sequence
$$\pi_0(G_i')\to \pi_0(G_i)\to \pi_0(H)=\{1\}$$
we get that $\pi_0(G_i')\to \pi_0(G_i)$ is onto. From this and diagram (\ref{Diag:connect}) we deduce that the map
$$\pi_0(K_i)\to \pi_0(G_i)$$
is onto. But $K_i$ and $G_i$ have the same number of connected components, namely 1,2,4 for $i=0,1,2$ respectively.
Hence this map is an isomorphism. \qed

\hfill

Let $\fg=\fk\oplus \fp$ be the Cartan decomposition of the complexified Lie algebra of all three groups $G_i$, $i=0,1,2$.
Explicitly we choose them to be
\begin{eqnarray*}
\fk=\left\{\left[\begin{array}{cc}
            X&0\\
            0&Y
            \end{array}\right]\big|\, X\in Mat_{p\times p}(\CC), Y\in Mat_{q\times q}(\CC)\mbox{ are anti-symmetric}\right\},\\
\fp=\left\{\left[\begin{array}{cc}
                  0_{p\times p}&iW\\
                  -iW^t&0_{q\times q}
                  \end{array}\right]\big|\, W\in Mat_{p\times q}(\CC)\right\}.
\end{eqnarray*}
Let us choose $\fa\subset \fp$ to be a maximal abelian subalgebra as follows:
\begin{eqnarray}\label{Def-a}
\fa=\left\{\left[\begin{array}{c|c|c}
                  0_{p\times p}&D&0\\\hline
                  -D&0_{p\times p}&0\\\hline
                  0&0&0_{q-p\times q-p}
                 \end{array}\right]\big|\, D\in Mat_{p\times p}(\CC) \mbox{ is complex diagonal}\right\}.
\end{eqnarray}

Let us choose the basis $\lam_1,\dots,\lam_p$ in $\fa^*$ as follows:
the value of $\lam_j$ on the element of the above form with
$D=i\cdot diag(a_1,\dots,a_p)$ is equal to $a_j$.

Then by \cite{loos}, Ch. VII, \S 2.3, the (non-zero) roots $\Sigma$
of $\fa$ in $\fg$ are
\begin{align}\label{Sigma1}
&\pm(\lam_i\pm\lam_j), \, 1\leq i<j\leq p, \mbox{ with multiplicity } 1;\\\label{Sigma2}
&\pm\lam_j, 1\leq j\leq p, \mbox{ with multiplicity } q-p.
\end{align}
Let us choose the positive roots as follows
\begin{eqnarray}\label{E:positive-roots}
\Sigma^+=\{-(\lam_i\pm\lam_j)|\, 1\leq i<j\leq p\}\cup \{-\lam_j|\,1\leq j\leq p\}.
\end{eqnarray}
The multiplicity of the root $\alp$ will be denoted by $m_\alp$.

The half sum $\rho$ of positive roots (counting multiplicities) is equal to
\begin{eqnarray}\label{E:half-sum}
\rho=-\sum_{i=1}^p(\frac{n}{2}-i)\lam_i.
\end{eqnarray}

Finally let us define $\fn$ to be the $\CC$-span of positive root spaces; $\fn\subset\fg$ is a nilpotent subalgebra.
Then
$$\fg=\fk\oplus \fa\oplus \fn$$
is the Iwasawa decomposition of the Lie algebra. Let $A,N\subset
G_0$ be the real analytic subgroups obtained by the exponentiation
of $\fa\cap so(p,q),\fn\cap so(p,q)$ respectively.
\begin{proposition}\label{P:1.5}
For $i=0,1,2$ we have the Iwasawa decomposition for groups $G_i$, namely the product map
$$K_i\times A\times N\to G_i$$
is a diffeomorphism.
\end{proposition}
{\bf Proof.} For $i=0$ this is a classical result since the group
$G_0$ is connected, see e.g. \cite{gga}, Ch. VI, Theorem 5.1,. In
general the result follows immediately from that case and
Proposition \ref{AP:1}. \qed

\hfill

Let us define
\begin{eqnarray*}
M_i:=\{g\in K_i|\, Ad(g)(a)=a\mbox{ for any } a\in \fa\},\\
M_i'=\{g\in K_i|\, Ad(g)(\fa)\subset \fa\}.
\end{eqnarray*}
Clearly $M_i\lhd M_i'$. The group $W_i:=M_i'/M_i$ is finite, as will be seen below, and it is called the little Weyl group
of the symmetric space $G_i/K_i$. Naturally $W_i\subset GL(\fa)$.
Let us describe $M_i, M_i'$ explicitly; we leave the details of the elementary computations to the reader:
\begin{eqnarray*}
M_2=\left\{\left[\begin{array}{c|c|c}
                   \eps&0&0\\\hline
                   0&\eps&0\\\hline
                   0&0&V
                  \end{array}\right]\right\},
M_2'=\left\{\left[\begin{array}{c|c|c}
                   \eps S&0&0\\\hline
                   0&\delta S&0\\\hline
                   0&0&V
                  \end{array}\right]\right\}
\end{eqnarray*}
where $S$ runs over all permutations of the standard basis in
$\RR^p$, $\eps,\delta$ run over diagonal  $p\times p$  matrices with
$ \pm 1 $ on the diagonal, $V\in O(q-p)$. Notice that if $p=q$ then
$V$ disappears from both formulas. The matrix
$\left[\begin{array}{c|c|c}
                   \eps S&0&0\\\hline
                   0&\delta S&0\\\hline
                   0&0&V
                  \end{array}\right]$ acts on an element $i\cdot diag(x_1,\dots,x_p)\in \fa$ by
\begin{eqnarray}\label{E:action}
i\cdot diag(x_1,\dots,x_p)\mapsto  (\eps\delta^{-1})\cdot i\cdot diag(x_{S(1)},\dots,x_{S(p)}).
\end{eqnarray}
From this the above description of $M_2$ is obvious.

It is clear that
$$M_i=M_2\cap K_i,\, M_i'=M_2'\cap K_i \mbox{ for } i=0,1.$$
Hence, in particular, $M_i'/M_0'\subset K_i/K_0$. But using the above descriptions of $M_i'$ it is easy to see that this
inclusion is in fact an equality:
\begin{eqnarray}\label{E:M-K}
M_i'/M_0'= K_i/K_0.
\end{eqnarray}

The action of $M_2'$ on $\fa$ induces a group homomorphism $M_2'\to GL(\fa)$ with kernel $M_2$. Restricting it to
$M_i'$, $i=0,1$ we get imbeddings
$$W_0\subset W_1\subset W_2.$$
Using these formulas it is easy to show that
\begin{align}\label{E:Weyl}
&W_0=W_1=W_2=\{\pm 1\}^p\rtimes S_p \mbox{ for } p<q,\\
&W_0=W_1=\{\pm 1\}^p_{even}\rtimes S_p,\, W_2=\{\pm 1\}^p\rtimes S_p \mbox{ for } p=q,
\end{align}
where $\{\pm 1\}^p_{even}$ denotes the index 2 subgroup of $\{\pm 1\}^p$ consisting sequences of $\pm 1$ with even number of minuses.

%$\rho$, $\Sigma^+$, $m_\alp$ were not defined !!!!!!!! *******************

\hfill

\begin{proposition}\label{P:2}
Let $K_i\cdot A\cdot N$ be the Iwasawa decomposition. Let $dk,da,dn$
be Haar measures on $K_i,A,N$ respectively (all of them are both
left and right invariant). Then the Haar measure $dg$ on $G_i$ can be
normalized so that
$$dg=e^{2\rho(\log \bullet)}dk\cdot da \cdot dn,$$
where $\log\colon A\to \fa\cap so(p,q)$ is the inverse of the
exponential map.
\end{proposition}
{\bf Proof.} For $G_0$ this is proven in Proposition 5.1, Ch. I, \S 5 in \cite{gga}. The cases $i=1,2$ follow from this one by applying Proposition \ref{AP:1}
since $G_0$ is the connected component of the identity of $G_i$, and $K_0$ of $K_i$. \qed

\begin{proposition}[Harish-Chandra]\label{P:3}
For $f\in C_c(G_i)$ we have
$$\int_{G_i}f(g)dg=\int_{K_iNA}f(kna)dk\cdot dn\cdot da=\int_{ANK_i}f(ank)da\cdot dn \cdot dk.$$
\end{proposition}
{\bf Proof.} For $i=0$ this is Corollary 5.3, Ch. I, \S 5, in
\cite{gga}. Combined with Proposition \ref{AP:1} it implies other
cases. \qed

\begin{proposition}\label{P:4}
For $H\in \fa\cap so(p,q)$, $a=\exp(H)$, let
$$D(a)=\prod_{\alp\in \Sigma^+}(\sinh(\frac{1}{2}\alp(H)))^{m_\alp}.$$
Then for a suitable normalization of invariant measure $dg_A$ on $G/A$ we have for any $a\in A$ such that $D(a)\ne 0$:
$$|D(a)|\cdot\int_{G/A}f(gag^{-1})dg_A=e^{\rho(\log a)}\int_{K\times N}f(kank^{-1})dk\cdot dn$$
for $f\in C_c(G_i)$ (the integrals on both sides are claimed to be absolutely convergent).
\end{proposition}
{\bf Proof.} For $G_0$ this is proven in Proposition 5.6 in Ch. I, \S 5 of \cite{gga}. The cases $i=1,2$ follows from this one and Proposition \ref{AP:1} by additional summation over connected components
of $G_i$, or equivalently $K_i$. \qed

\begin{proposition}\label{AP:5}
Let $f\in C_c(G_i)$ satisfy $f(kgk^{-1})=f(g)$ for any $k\in K_i,g\in G_i$.  Then the function
$$F_f(a):=e^{\rho(\log a)}\int_Nf(an)dn,\, a\in A,$$
satisfies $$F_f(a^s)=F_f(a),\, \mbox{ for each } a\in A, s\in W_i.$$
\end{proposition}
{\bf Proof.} We just repeat the argument of Proposition 5.7 in Ch. I, \S 5 in \cite{gga}. By continuity it suffices to prove the identity for $D(a)\ne 0$. By
Proposition \ref{P:4} we have
$$F_f(a)=|D(a)|\int_{G/A}f(gag^{-1})dg_A.$$
We have to show that the right hand side of the last equality is $W_i$-invariant; we will do it in fact for any $f\in C_c(G_i)$.

$W_i$ permutes the (restricted) roots $\Sigma$. Hence $W_i$
preserves $|D(a)|$. Let $u\in M_i'$ be a representative of $s$.
Since $uAu^{-1}=A$, the map $\phi\colon G_i/A\to G_i/A$ given by
$\phi(xA)=uxu^{-1}A$ is a well defined diffeomorphism. The
conjugation  by $u$ preserves Haar measures on $G_i$ and on $A$
since $u$ is contained in the finite group $W_i$, hence $\phi$
preserves the $G_i$-invariant measure on $G_i/A$. We have, using $a^s=uau^{-1}$,
\begin{align*}
&\int_{G_i/A}f(ga^sg^{-1})dg_A=\int_{G_i/A}f((ugu^{-1})a^s(ugu^{-1})^{-1})dg_A\\
&=\int_{G_i/A}f(ugag^{-1}u^{-1})\overset{\mbox{Prop. } \ref{P:3}}{=}\int_{K_i\times N}f(uknan^{-1}k^{-1}u^{-1})dk\cdot dn\\
&=\int_{K_i\times N}f(knan^{-1}k^{-1})dk\cdot dn=\int_{G/A}f(gag^{-1}) dg_A.
\end{align*}
The proposition is proved. \qed

\hfill

For $D\in U(\fg)\simeq D(G_i)$ let us denote by $D_\fa\in U(\fa)=D(A)$ the image of $D$ under the projection
$$\tilde\gamma\colon U(\fg)=U(\fn)\otimes U(\fa)\otimes U(\fk)\to U(\fa)$$
which is just the quotient map under the subspace $\fn\cdot U(\fg)+U(\fg)\cdot \fk$.
By Lemma 5.14 in Ch. II in \cite{gga} for any $\phi\in C_c(G_i)$ such that
$$\phi(ngk)=\phi(g) \mbox{ for any } g\in G_i, n\in N,k\in K_i$$
one has
\begin{eqnarray}\label{Eq:1}
(D\phi)|_A=D_{\fa}(\phi|_A).
\end{eqnarray}
(Strictly speaking, in \cite{gga} the equality (\ref{Eq:1}) is stated only for connected groups, i.e. for $G_0$.
But to prove (\ref{Eq:1}) for $G_i$ it suffices to restrict $\phi$ to $G_0$.)

For $g\in G_i$ we denote by $A(g)\in \fa\cap so(p,q)$ the unique
element such that $g\in N\cdot \exp(A(g))\cdot K_i$.
\begin{lemma}\label{L:6}
For each linear functional $\nu\colon \fa\to \CC$ the function $$\phi(g):=\int_{K_i}e^{\nu(A(kg))}dk$$ satisfies
for any $D\in D_{K_i}(G_i)=U(\fg)^{K_i}$
$$D\phi=D_\fa(\nu)\phi,$$
where $D_\fa(\nu)\in \CC$ is defined as follows: we have
$$D_\fa\in D(\fa)=U(\fa)=S(\fa)=\CC[\fa^*],$$
where $\CC[\fa^*]$ denotes complex polynomials on $\fa^*$; thus $D_\fa$ can be evaluated at $\nu\in \fa^*$.
\end{lemma}
{\bf Proof.} Let $F(g):=e^{\nu(A(g))}$. Clearly $F(ngk)=F(g)$. Hence by (\ref{Eq:1}) for any $a\in A$
\begin{eqnarray}\label{Eq:2}
(DF)(a)=(D_\fa F|_A)(a)=D_\fa(\nu)\cdot F(a).
\end{eqnarray}
The functions $DF$ and $D_\fa(\nu)\cdot F$ are both left $N$-invariant and right $K_i$-invariant, and by (\ref{Eq:2}) coincide on $A$. Hence
$$DF=D_\fa(\nu)\cdot F.$$
Hence
\begin{eqnarray*}
(D\phi)(g)=D_g\left(\int_{K_i}F(kg)dk\right)=\int_{K_i}(DF)(kg)dk=\\
D_\fa(\nu)\int_{K_i}F(kg)dk=D_\fa(\nu)\cdot \phi(g).
\end{eqnarray*}
\qed

\begin{corollary}\label{C:6.5}
The map $D_{K_i}(G_i)= U(\fg)^{K_i} \to  D(\fa)\simeq U(\fa)$
given by
$$D\mapsto D_\fa$$
is a homomorphism of algebras.
\end{corollary}
{\bf Proof.} In the notation of the proof of Lemma \ref{L:6} for any $\nu\in \fa^*$ the corresponding
function $\phi$ is an eigenfunction of any $D\in D_{K_i}(G_i)$ with the eigenvalue $D_\fa(\nu)$. Hence for any $D_1,D_2\in D_{K_i}(G_i)$
this implies that
$$(D_1\circ D_2)_\fa(\nu)=D_{1\fa}(\nu)\cdot D_{2\fa}(\nu).$$
Since this holds for any $\nu$ we get $(D_1\circ D_2)_\fa=D_{1\fa}\circ D_{2\fa}$. The corollary follows. \qed

\begin{theorem}\label{T:7}
With $\rho=\frac{1}{2}\sum_{\alp\in \Sigma^+}m_\alp\cdot \alp$ as above, for $\nu\in \fa^*$ let us denote
$$\phi_\nu(g):=\int_{K_i}e^{(\nu+\rho)(A(kg))}dk, \, g\in G_i.$$
Then $\phi_{s\nu}=\phi_\nu$ for each $s\in W_i$, where $(s\nu)(H)=\nu(s^{-1}H),\, H\in \fa$.
\end{theorem}
{\bf Proof.} Since $\phi_\nu, \phi_{s\nu}$ are bi-invariant under $K_i$, it suffices to prove for any $f\in C_c(G_i)$
which is bi-invariant under $K_i$ that
\begin{eqnarray}\label{Eq:3}
\int_{G_i}\phi_{s\nu}(g)f(g)dg=\int_{G_i}\phi_\nu(g)f(g)dg.
\end{eqnarray}
Under the decomposition $G_i=A\cdot N\cdot K_i$ by Proposition \ref{P:3} we have
\begin{eqnarray*}
dg=da\cdot dn\cdot dk.
\end{eqnarray*}
Hence for $f\in C_c(G_i)$ being $K_i$-bi-invariant we have
\begin{eqnarray*}
\int_{G_i}\phi_\nu(g)f(g)dg&=\int_{K_i}dk\int_{G_i}e^{(\nu+\rho)(A(kg))}f(g)dg\\
=\int_{G_i}e^{(\nu+\rho)(A(g))}f(g)dg&=\int_A da\int_N e^{(\nu+\rho)(\log a)}f(an)dn\\
&=\int_Ae^{\nu(\log a)}F_f(a)da,
\end{eqnarray*}
where $F_f(a)=e^{\rho(\log a)}\int_Nf(an) dn$. But by Proposition \ref{AP:5} $F_f(a^s)=F_f(a)$. Hence
\begin{eqnarray*}
\int_{G_i}\phi_{s\nu}(g)f(g)dg=\int_A e^{(s\nu)(\log a)}F_f(a)da=\\
\int_Ae^{\nu(\log a)}F_f(a^s)da=\int_{G_i}\phi_\nu(g)f(g)da.
\end{eqnarray*}
\qed

\hfill

Let us define
\begin{eqnarray}\label{Def:gamma}
\gamma\colon D_{K_i}(G_i)=U(\fg)^{K_i}\to U(\fa)
\end{eqnarray}
by $\gamma(D)=e^{-\rho}\circ D_\fa\circ e^\rho$, where $e^\rho$ denotes the operator of shift in $A$ by the element $e^\rho$.
From Lemma \ref{L:6}, Corollary \ref{C:6.5}, Theorem \ref{T:7}, and Proposition \ref{P:homomor-operators} we immediately get
\begin{corollary}\label{AC:8}
$$\gamma\colon D_{K_i}(G_i)\to U(\fa)$$
is a homomorphism of algebras. Moreover its image is contained in $U(\fa)^{W_i}$, and the kernel is equal to
$U(\fg)^{K_i}\cap (U(\fg)\cdot \fk)$.
\end{corollary}
\begin{remark}\label{R:remark-Weyl}
More explicitly, for any $D\in D_{K_i}(G_i)$, any $\nu\in \fa^*$, and any $s\in W_i$ one has
$$D_\fa(s\nu+\rho)=D_\fa(\nu+\rho),$$
and $\gamma(D)(\nu)=D_\fa(\nu+\rho)$.
\end{remark}

Our first main goal is to prove
\begin{theorem}\label{T:HC}
For $i=0,1,2$
$$\gamma\colon D_{K_i}(G_i)=U(\fg)^{K_i}\to U(\fa)^{W_i}$$
is an epimorphism of algebras with the kernel $U(\fg)^{K_i}\cap (U(\fg)\cdot \fk)$.
\end{theorem}
It remains only to prove that $\gamma$ is onto.
We will need a lemma. Let us denote for brevity
$$\fa_0:=\fa\cap so(p,q),\, \fp_0:=\fp\cap so(p,q).$$
\begin{lemma}\label{AL:9}
The restriction map
$$C^\infty(\fp_0)^{K_i}\to C^\infty(\fa_0)^{W_i}$$
is an isomorphism of algebras.
\end{lemma}
{\bf Proof.} For $i=0$ this is Corollary 5.11(i) in Ch. II of \cite{gga}. For $i=1,2$ the injectivity follows from the case $i=0$.
Hence it remains to prove surjectivity.

Let $f\in C^\infty(\fa_0)^{W_i}\subset C^\infty(\fa_0)^{W_0}$. Then by case $i=0$ there exists $\tilde F\in C^\infty(\fp_0)^{K_0}$ such that
$$\tilde F|_{\fa_0}=f.$$
Let us define a new function on $\fp_0$ (below the Haar measure $dk$ on $K_i$ is normalized to be probability measure):
$$F(x)=\int_{K_i}\tilde F(k^{-1}xk)dk=\frac{1}{|K_i/K_0|}\sum_{\sigma\in K_i/K_0} \tilde F(\sigma^{-1}x \sigma).$$
Clearly $F\in C^\infty(\fp_0)^{K_i}$. Let us show that $F|_{\fa_0} =f$.

By (\ref{E:M-K}) one has $K_i/K_0=M_i'/M_0'$. Then for any $a\in \fa_0$ we have
\begin{eqnarray*}
F(a)=\frac{1}{|M_i'/M_0'|}\sum_{\sigma\in M_i'/M_0'} \tilde F(\sigma^{-1}a \sigma)=f(a),
\end{eqnarray*}
where the last equality is due to the facts that $\sigma^{-1}a \sigma\in \fa_0$ and $f$ is $W_i$-invariant. \qed

\begin{corollary}\label{C:13}
The restriction of polynomials
$$\CC[\fp]^{K_i}\to \CC[\fa]^{W_i}$$
is an isomorphism of algebras.
\end{corollary}
{\bf Proof.} By Lemma \ref{AL:9} only surjectivity has to be proven.
Let $f\in \CC[\fa]^{W_i}$ be a homogeneous polynomial. By Lemma \ref{AL:9} there exists
unique $F\in C^\infty[\fp_0]^{K_i}$ such that $F|_{\fa_0}=f$. $F$ is also homogeneous of the same degree as $f$.
But any infinitely smooth homogeneous function is a polynomial, i.e. $F\in \CC[\fp]^{K_i}$. \qed

\hfill

{\bf Proof of Theorem \ref{T:HC}.} Let us define a linear map
$$\lam\colon Sym(\fg)\to U(\fg)$$
by $\lam(X_1\otimes\dots\otimes X_m):=\frac{1}{m!}\sum_{\sigma\in S_m}X_{\sigma(1)}\dots X_{\sigma(m)}\in U(\fg)$.
It is easy to see that $\lam$ is an isomorphism of vector spaces which commutes with the adjoint action of $G_i$.

Let $\fq$ be the orthogonal complement of $\fa$ in $\fp$ with respect to the Killing form (recall that the restriction of the Killing form to
$\fp_0$ is positive definite). Thus
$$\fp=\fa\oplus \fq.$$
Since the Killing form on $\fg$ is $G_i$-invariant, $\fq$ is $M_i'$-invariant. The Killing form induces identifications
$$\fp^*\simeq \fp,\, \fa^*\simeq \fa$$
such that the first is $K_i$-equivariant, and the second is $M_i'$-equivariant. Under these identifications the restriction map
$$\CC[\fp]\to \CC[\fa]$$
becomes the projection
$$\tau\colon Sym(\fp)=Sym(\fa)\oplus Sym(\fp)\cdot \fq\to Sym(\fa).$$

By \cite{gga}, Ch. II, formula (38), for any $p\in Sym(\fp)^{K_0}$ one has
\begin{eqnarray}\label{E:10}
\mbox{degree}(\gamma(\lam(p))-\tau(p))<\mbox{degree}(p).
\end{eqnarray}

Let us fix now $t\in U(\fa)^{W_i}=Sym(\fa)^{W_i}$. We have to show that $t$ belongs to the image of the homomorphism $\gamma$.
We may assume that $t$ is homogeneous. By Corollary \ref{C:13} and the above remarks on equivariant identifications, there exists unique $p\in Sym(\fp)^{K_i}$
such that $\tau(p)=t$. Also $p$ has the same degree as $t$. By (\ref{E:10})
\begin{eqnarray}\label{E:11}
\mbox{degree}(\gamma(\lam(p))-t)<\mbox{degree}(p).
\end{eqnarray}

But $\lam(p)\in U(\fg)^{K_i}$. Continuing by induction in $\mbox{degree}(t)$ we prove the theorem. \qed

\hfill

From Theorem \ref{T:HC} and Proposition \ref{P:homomor-operators} we immediately deduce
\begin{corollary}\label{C:14}
The epimorphism $\gamma$ induces the isomorphism, also denoted by $\gamma$ and called the Harish-Chandra isomorphism,
$$\gamma\colon D(G_i/K_i)\tilde\to U(\fa)^{W_i}.$$
\end{corollary}

\hfill

Now let us start discussing the differential operators on compact groups. For $n=p+q$ let us denote
$$U_0=U_1=SO(n),\, U_2=O(n).$$
Then $K_i\subset U_i$ for $i=0,1,2$. Moreover $U_0/K_0$ is the Grassmannian of oriented $p$-planes in $\RR^n$, while $U_1/K_1=U_2/K_2$
is the Grassmannian of $p$-planes in $\RR^n$. Notice that the complexified Lie algebra of $U_i$ is naturally identified with $\fg=so(n,\CC)$.
By Proposition \ref{P:homomor-operators} the algebra $D(U_i/K_i)$ of $U_i$-invariant operators on $U_i/K_i$ is naturally identified with quotient algebra $U(\fg)^{K_i}/(U(\fg)^{K_i}\cap (U\left(\fg)\cdot\fk)\right)$.
Hence we obtain an isomorphism
\begin{eqnarray}\label{E:12}
D(U_i/K_i)\simeq D(G_i/K_i).
\end{eqnarray}

In order to finish the proof of Theorem \ref{HC} it remains to prove the following proposition.
\begin{proposition}\label{P:compact-repr}
Let $n=p+q,\, r=\min\{p,q\}$. Let $\cv\subset L^2(O(n)/O(p)\times O(q))$ be an irreducible representation corresponding to a sequence
$m$ of even non-negative integers
$$m=(m_1\geq \dots\geq m_{r-1}\geq m_r\geq 0).$$
Let $D$ be an invariant operator on $O(n)/O(p)\times O(q)$. Then $D$ acts on $\cv$ by multiplication by a scalar $\gamma(D)(m+\bar\rho)$, where $\bar\rho$ is the sequence given by
\begin{eqnarray}\label{E:bar-rho}
\bar\rho_i=n/2-i,\, i=1,\dots,r.
\end{eqnarray}
\end{proposition}

Before the proof let us introduce some more notation. Let us agree as previously $0<p\leq q$. Let $\fh\subset \fg$ be a $\theta$-invariant Cartan subalgebra containing $\fa$; we do not need to specify it more explicitly.
Thus
$$\fh=\fa\oplus \fb,$$
where $\fb=\fh\cap \fk$. Moreover $\fh$ can be chosen to be of the form
$$\fh=\fh_0\otimes_\RR\CC,$$
where $\fh_0=\fa_0\oplus \fb_0$ and $\fa_0=\fa\cap so(p,q),\, \fb_0=i(\fb\cap (so(p)\times so(q))).$ Thus in any unitary representation of $U_i$ the algebra $\fh_0$ acts by
self-adjoint operators. In particular the restriction to $\fh_0$ of any root of $\fg$ with respect to $\fh$ is a real valued linear functional on $\fh_0$.

Let us choose the positive roots system of $\fg$ with respect to $\fh$ such that the restriction of these roots to $\fa$ is equal to $\Sigma\backslash \Sigma^+=-\Sigma^+$, where $\Sigma$ is defined by
(\ref{Sigma1})-(\ref{Sigma2}), and $\Sigma^+$ is defined by (\ref{E:positive-roots}).
Let $\fn^+\subset\fg$ be the $\CC$-span of these roots. Thus $\fh\oplus \fn^+$ is a Borel subalgebra of $\fg$. Then we have
\begin{eqnarray}\label{E:nilp-inclusion}
\theta(\fn)\subset \fn^+.
\end{eqnarray}
Indeed for $\alp\in \Sigma$ if
$$[a,x]=\alp(a)x \mbox{ for any } a\in\fa,$$
then applying the automorphism $\theta$ on both sides and using $\theta(a)=-a$ we get
$$[a,\theta(x)]=-\alp(a)\theta(x).$$
This implies (\ref{E:nilp-inclusion}).

\hfill

%Part (1) of the following lemma is an exact analogue for disconnected groups of Lemma 3 on p. 92 in \cite{takeuchi}.
\begin{lemma}\label{L:non-orthog}
Let $\cv$ be an irreducible representation of the group $SO(n)$. Let $v\in\cv$ be a highest weight vector, i.e. a non-zero vector satisfying
\begin{eqnarray}\label{E:hw-vector}
(\fh\oplus \fn^+)(v)\subset \CC\cdot v.
\end{eqnarray}
%(Note that since the group $O(n)$ is not connected the vector $V$ may not be unique (up to a constant).)
Let $w\in \cv$ be a non-zero $S(O(p)\times O(q))$-invariant vector.
Let $(\cdot,\cdot)$ be a $SO(n)$-invariant hermitian product. Then

(1) $(v,w)\ne 0.$

(2) $\fb (v)=0$.
\end{lemma}
{\bf Proof.} %Let us assume in the contrary that $(V,W)=0$.
%The representation of $O(n)$ extends uniquely to an algebraic representation of the complex group $O(n,\CC)$.
%Then we can restrict it to $O(p,q)$. This restriction is still irreducible since $O(p,q)$ and $O(n)$ have the same complexification $O(n,\CC)$.
%Consider the Iwasawa decomposition $O(p,q)=(O(p)\times O(q))\cdot A\cdot N$. Then by our assumption
%$$A\cdot N(V)\subset \CC\cdot V.$$
%Hence $O(p,q)(V)\subset \CC \cdot (O(p)\times O(q))(V).$ But if $(V,W)=0$ then $(O(p,q))(V)\subset W^\perp$. Hence subspace generated by $V$ is contained in $W^\perp$, in particular is proper invariant.
%This is a contradiction. Hence part (1) is proved.

Part (1) is a special case of Lemma 3 on p. 92 in \cite{takeuchi}.

To prove part (2) it suffices to show that for any $X\in \fb_0$ one has $X(v)=0$. But by (\ref{E:hw-vector})
$$X(v)=c\cdot v$$
for some $c\in \CC$. Since $X$ is a self-adjoint operator
$$c(v,w)=(X(v),w)=(v,X(w))=0,$$
where the last equality is by assumption that $w$ is $S(O(p)\times O(q))$-invariant. Now part (1) implies part (2). \qed

{\bf Proof of Proposition \ref{P:compact-repr}.} Let us restrict our representation $\cv$ to $SO(n)$. If either $r<n/2$, or $r=n/2$ and $m_r=0$ then $\cv$ is an irreducible $SO(n)$-module.
Otherwise $\cv$ is a direct sum of two irreducible $SO(n)$-modules with
highest weights $(m_1,\dots,m_{\frac{n}{2}-1},m_{\frac{n}{2}})$ and $(m_1,\dots,m_{\frac{n}{2}-1},-m_{\frac{n}{2}})$. Hence it suffices to prove that for any irreducible $SO(n)$-module
$$\cv_1\subset L^2(SO(n)/S(O(p)\times O(q))$$ with the highest weight
$$m_1\geq\dots\geq m_{\frac{n}{2}-1}\geq m_{\frac{n}{2}}\geq 0$$
the action of $SO(n)$-invariant operator $D$ on $\cv_1$ is multiplication by the scalar $\gamma(D)(m+\bar\rho)$, where $\bar\rho$ is given by (\ref{E:bar-rho}).

Let us choose $L\in U(\fg)^{S(O(p)\times O(q))}$ to be a representative of $D$. Let $(\cdot,\cdot)$ be an $SO(n)$-invariant hermitian product which is assumed to be linear with respect to the second argument and
anti-linear with respect to the first.

Let $v\in \cv_1$ be a highest weight vector, i.e.
\begin{eqnarray}\label{E:hw}
\fn^+v=0,\, \fh v\subset \CC\cdot v.
\end{eqnarray}
Let $\lam\in \fh^*$ be the corresponding highest weight, i.e.
$$h(v)=\lam(h) v\,\, \forall h\in \fh.$$
Note that by Lemma \ref{L:non-orthog}(2), $\lam$ vanishes on $\fb$, i.e.
\begin{eqnarray}\label{E:lambda-hw}
\lam\in \fa^*\subset \fh^*.
\end{eqnarray}

Let $w\in \cv_1$ be an $S(O(p)\times O(q))$-invariant non-zero vector. By Lemma \ref{L:non-orthog}(1)
\begin{eqnarray}\label{E:non-zero}
(v,w)\ne 0.
\end{eqnarray}
Hence it suffices to show that
$$(v,Dw)=(v,w)\cdot \gamma(D)(m+\bar\rho).$$

It is clear that for any $D_1\in U(\fg)\cdot \fk$ one has $D_1 w=0$ hence
\begin{eqnarray}\label{E:vanish1}
(v,D_1(w))=0.
\end{eqnarray}
Next we claim that for $D_2\in \fn\cdot U(\fg)$ one has
\begin{eqnarray}\label{E:adj-vanish}
D_2^*v=0.
\end{eqnarray}
Hence it follows that
\begin{eqnarray}\label{E:vanish2}
(v,D_2(w))=0.
\end{eqnarray}
Indeed, to prove (\ref{E:adj-vanish}) let us choose $n\in \fn$. We may assume that $n\in \fn\cap so(p,q)$. Then using the Cartan decomposition for $so(p,q)$ consider
$$n=n_1+n_2,$$
where $n_1\in \fk\cap so(p,q),n_2\in \fp\cap so(p,q)$.  The action on $n_1$ in $\cv_1$ is anti-symmetric, while the action of $n_2$ on $\cv_1$ is symmetric. Hence
$$n^*=-n_1+n_2=(-\theta)(n)\in \fn^+,$$
where the last inclusion is by (\ref{E:nilp-inclusion}).
Hence (\ref{E:adj-vanish}) follows.

Then (\ref{E:vanish1}) and (\ref{E:vanish2}) imply that
\begin{eqnarray*}
(v,Dw)=(v,D_{\fa}w)=((D_\fa)^*v,w)=\\
D_\fa (\lam)(v,w)=\gamma(D)(\lam-\rho)\cdot(v,w),
\end{eqnarray*}
where the last equality is due to Remark \ref{R:remark-Weyl}.
Since $(v,w)\ne 0$ we get that $D$ acts on $\cv_1$, and hence on $\cv$, by the multiplication by the scalar $\gamma(D)(\lam-\rho)$.
We remind that here $\lam\in \fa^*\subset\fh^*$, and $\rho\in \fa^*$ is the half sum of (restricted) roots from $\Sigma^+$ given by (\ref{E:positive-roots}).

It remains to translate the last result into the combinatorial language. First we identify $\fa^*\tilde\to \CC^p$ as follows. First identify $\fa$ with $\CC^p$ as in (\ref{Def-a})
with the matrix $D$ being diagonal with entries from $\CC$. Dualizing this isomorphism we get an identification $\fa^*\simeq \CC^p$.

If $\lam\simeq \fa^*\subset \fh^*$ is the highest weight of this representation with the above choice of $\fn^+$ then it corresponds to a sequence
\begin{eqnarray*}
m_1\geq\dots\geq m_{p-1}\geq |m_p| \mbox{ if } p=\frac{n}{2},\\
m_1\geq\dots\geq m_{p-1}\geq m_p\geq 0 \mbox{ if } p<\frac{n}{2}.
\end{eqnarray*}
Recall that in our case we have chosen $m_p\geq 0$ from the very beginning.

Furthermore under the above identification $\rho$ corresponds to the sequence $\tilde\rho$ with
$$\tilde\rho_i=-(\frac{n}{2}-i), \, i=1,\dots,p.$$
Clearly $\tilde \rho=-\bar\rho$, where $\bar \rho$ is defined by (\ref{E:bar-rho}).
Theorem is proved.
\qed

\end{document}

%% file: diagram.tex
\makeatletter

\def\diagram{\m@th\leftwidth=\z@ \rightwidth=\z@ \topheight=\z@
\botheight=\z@ \setbox\@picbox\hbox\bgroup}

\def\enddiagram{\egroup\wd\@picbox\rightwidth\unitlength
\ht\@picbox\topheight\unitlength \dp\@picbox\botheight\unitlength
\hskip\leftwidth\unitlength\box\@picbox}

\def\bfig{\begin{diagram}}
\def\efig{\end{diagram}}
\newcount\wideness \newcount\leftwidth \newcount\rightwidth
\newcount\highness \newcount\topheight \newcount\botheight

\def\ratchet#1#2{\ifnum#1<#2 \global #1=#2 \fi}

\def\putbox(#1,#2)#3{%
\horsize{\wideness}{#3} \divide\wideness by 2
{\advance\wideness by #1 \ratchet{\rightwidth}{\wideness}}
{\advance\wideness by -#1 \ratchet{\leftwidth}{\wideness}}
\vertsize{\highness}{#3} \divide\highness by 2
{\advance\highness by #2 \ratchet{\topheight}{\highness}}
{\advance\highness by -#2 \ratchet{\botheight}{\highness}}
\put(#1,#2){\makebox(0,0){$#3$}}}

\def\putlbox(#1,#2)#3{%
\horsize{\wideness}{#3}
{\advance\wideness by #1 \ratchet{\rightwidth}{\wideness}}
{\ratchet{\leftwidth}{-#1}}
\vertsize{\highness}{#3} \divide\highness by 2
{\advance\highness by #2 \ratchet{\topheight}{\highness}}
{\advance\highness by -#2 \ratchet{\botheight}{\highness}}
\put(#1,#2){\makebox(0,0)[l]{$#3$}}}

\def\putrbox(#1,#2)#3{%
\horsize{\wideness}{#3}
{\ratchet{\rightwidth}{#1}}
{\advance\wideness by -#1 \ratchet{\leftwidth}{\wideness}}
\vertsize{\highness}{#3} \divide\highness by 2
{\advance\highness by #2 \ratchet{\topheight}{\highness}}
{\advance\highness by -#2 \ratchet{\botheight}{\highness}}
\put(#1,#2){\makebox(0,0)[r]{$#3$}}}

\def\adjust[#1]{} % For compatibility

\newcount \coefa
\newcount \coefb
\newcount \coefc
\newcount\tempcounta
\newcount\tempcountb
\newcount\tempcountc
\newcount\tempcountd
\newcount\xext
\newcount\yext
\newcount\xoff
\newcount\yoff
\newcount\gap%
\newcount\arrowtypea
\newcount\arrowtypeb
\newcount\arrowtypec
\newcount\arrowtyped
\newcount\arrowtypee
\newcount\height
\newcount\width
\newcount\xpos
\newcount\ypos
\newcount\run
\newcount\rise
\newcount\arrowlength
\newcount\halflength
\newcount\arrowtype
\newdimen\tempdimen
\newdimen\xlen
\newdimen\ylen
\newsavebox{\tempboxa}%
\newsavebox{\tempboxb}%
\newsavebox{\tempboxc}%

\newdimen\w@dth

\def\setw@dth#1#2{\setbox\z@\hbox{\m@th$#1$}\w@dth=\wd\z@
\setbox\@ne\hbox{\m@th$#2$}\ifnum\w@dth<\wd\@ne \w@dth=\wd\@ne \fi
\advance\w@dth by 1.2em}

\def\t@^#1_#2{\allowbreak\def\n@one{#1}\def\n@two{#2}\mathrel
{\setw@dth{#1}{#2}
\mathop{\hbox to \w@dth{\rightarrowfill}}\limits
\ifx\n@one\empty\else ^{\box\z@}\fi
\ifx\n@two\empty\else _{\box\@ne}\fi}}
\def\t@@^#1{\@ifnextchar_{\t@^{#1}}{\t@^{#1}_{}}}
\def\to{\@ifnextchar^{\t@@}{\t@@^{}}}

\def\t@left^#1_#2{\def\n@one{#1}\def\n@two{#2}\mathrel{\setw@dth{#1}{#2}
\mathop{\hbox to \w@dth{\leftarrowfill}}\limits
\ifx\n@one\empty\else ^{\box\z@}\fi
\ifx\n@two\empty\else _{\box\@ne}\fi}}
\def\t@@left^#1{\@ifnextchar_{\t@left^{#1}}{\t@left^{#1}_{}}}
\def\toleft{\@ifnextchar^{\t@@left}{\t@@left^{}}}

\def\two@^#1_#2{\allowbreak
\def\n@one{#1}\def\n@two{#2}\mathrel{\setw@dth{#1}{#2}
\mathop{\vcenter{\lineskip\z@\baselineskip\z@
                 \hbox to \w@dth{\rightarrowfill}%
                 \hbox to \w@dth{\rightarrowfill}}%
       }\limits
\ifx\n@one\empty\else ^{\box\z@}\fi
\ifx\n@two\empty\else _{\box\@ne}\fi}}
\def\tw@@^#1{\@ifnextchar _{\two@^{#1}}{\two@^{#1}_{}}}
\def\two{\@ifnextchar ^{\tw@@}{\tw@@^{}}}

\def\tofr@^#1_#2{\def\n@one{#1}\def\n@two{#2}\mathrel{\setw@dth{#1}{#2}
\mathop{\vcenter{\hbox to \w@dth{\rightarrowfill}\kern-1.7ex
                 \hbox to \w@dth{\leftarrowfill}}%
       }\limits
\ifx\n@one\empty\else ^{\box\z@}\fi
\ifx\n@two\empty\else _{\box\@ne}\fi}}
\def\t@fr@^#1{\@ifnextchar_ {\tofr@^{#1}}{\tofr@^{#1}_{}}}
\def\tofro{\@ifnextchar^ {\t@fr@}{\t@fr@^{}}}

\def\mon{\mathop{\m@th\hbox to
      14.6\P@{\lasyb\char'51\hskip-2.1\P@$\arrext$\hss
$\mathord\rightarrow$}}\limits} % width of \epi
\def\leftmono{\mathrel{\m@th\hbox to
14.6\P@{$\mathord\leftarrow$\hss$\arrext$\hskip-2.1\P@\lasyb\char'50%
}}\limits} % width of \epi
\mathchardef\arrext="0200       % amr minus for arrow extension (see \into)

\def\settypes(#1,#2,#3){\arrowtypea#1 \arrowtypeb#2 \arrowtypec#3}
\def\settoheight#1#2{\setbox\@tempboxa\hbox{#2}#1\ht\@tempboxa\relax}%
\def\settodepth#1#2{\setbox\@tempboxa\hbox{#2}#1\dp\@tempboxa\relax}%
\def\settokens`#1`#2`#3`#4`{%
     \def\tokena{#1}\def\tokenb{#2}\def\tokenc{#3}\def\tokend{#4}}
\def\setsqparms[#1`#2`#3`#4;#5`#6]{%
\arrowtypea #1
\arrowtypeb #2
\arrowtypec #3
\arrowtyped #4
\width #5
\height #6
}
\def\setpos(#1,#2){\xpos=#1 \ypos#2}

\def\settriparms[#1`#2`#3;#4]{\settripairparms[#1`#2`#3`1`1;#4]}%

\def\settripairparms[#1`#2`#3`#4`#5;#6]{%
\arrowtypea #1
\arrowtypeb #2
\arrowtypec #3
\arrowtyped #4
\arrowtypee #5
\width #6
\height #6
}

\def\resetparms{\settripairparms[1`1`1`1`1;500]\width 500}%default values%

\resetparms

\def\mvector(#1,#2)#3{%%
\put(0,0){\vector(#1,#2){#3}}%
\put(0,0){\vector(#1,#2){26}}%
}
\def\evector(#1,#2)#3{{%%
\arrowlength #3
\put(0,0){\vector(#1,#2){\arrowlength}}%
\advance \arrowlength by-30
\put(0,0){\vector(#1,#2){\arrowlength}}%
}}

\def\horsize#1#2{%
\settowidth{\tempdimen}{$#2$}%
#1=\tempdimen
\divide #1 by\unitlength
}

\def\vertsize#1#2{%
\settoheight{\tempdimen}{$#2$}%
#1=\tempdimen
\settodepth{\tempdimen}{$#2$}%
\advance #1 by\tempdimen
\divide #1 by\unitlength
}

\def\putvector(#1,#2)(#3,#4)#5#6{{%
\ifnum3<\arrowtype
\putdashvector(#1,#2)(#3,#4)#5\arrowtype
\else
\ifnum\arrowtype<-3
\putdashvector(#1,#2)(#3,#4)#5\arrowtype
\else
\xpos=#1
\ypos=#2
\run=#3
\rise=#4
\arrowlength=#5
\ifnum \arrowtype<0
    \ifnum \run=0
        \advance \ypos by-\arrowlength
    \else
        \tempcounta \arrowlength
        \multiply \tempcounta by\rise
        \divide \tempcounta by\run
        \ifnum\run>0
            \advance \xpos by\arrowlength
            \advance \ypos by\tempcounta
        \else
            \advance \xpos by-\arrowlength
            \advance \ypos by-\tempcounta
        \fi
    \fi
    \multiply \arrowtype by-1
    \multiply \rise by-1
    \multiply \run by-1
\fi
\ifcase \arrowtype
\or \put(\xpos,\ypos){\vector(\run,\rise){\arrowlength}}%
\or \put(\xpos,\ypos){\mvector(\run,\rise)\arrowlength}%
\or \put(\xpos,\ypos){\evector(\run,\rise){\arrowlength}}%
\fi\fi\fi
}}

\def\putsplitvector(#1,#2)#3#4{%%
\xpos #1
\ypos #2
\arrowtype #4
\halflength #3
\arrowlength #3
\gap 140
\advance \halflength by-\gap
\divide \halflength by2
\ifnum\arrowtype>0
   \ifcase \arrowtype
   \or \put(\xpos,\ypos){\line(0,-1){\halflength}}%
       \advance\ypos by-\halflength
       \advance\ypos by-\gap
       \put(\xpos,\ypos){\vector(0,-1){\halflength}}%
   \or \put(\xpos,\ypos){\line(0,-1)\halflength}%
       \put(\xpos,\ypos){\vector(0,-1)3}%
       \advance\ypos by-\halflength
       \advance\ypos by-\gap
       \put(\xpos,\ypos){\vector(0,-1){\halflength}}%
   \or \put(\xpos,\ypos){\line(0,-1)\halflength}%
       \advance\ypos by-\halflength
       \advance\ypos by-\gap
       \put(\xpos,\ypos){\evector(0,-1){\halflength}}%
   \fi
\else \arrowtype=-\arrowtype
   \ifcase\arrowtype
   \or \advance \ypos by-\arrowlength
       \put(\xpos,\ypos){\line(0,1){\halflength}}%
       \advance\ypos by\halflength
       \advance\ypos by\gap
       \put(\xpos,\ypos){\vector(0,1){\halflength}}%
   \or \advance \ypos by-\arrowlength
       \put(\xpos,\ypos){\line(0,1)\halflength}%
       \put(\xpos,\ypos){\vector(0,1)3}%
       \advance\ypos by\halflength
       \advance\ypos by\gap
       \put(\xpos,\ypos){\vector(0,1){\halflength}}%
   \or \advance \ypos by-\arrowlength
       \put(\xpos,\ypos){\line(0,1)\halflength}%
       \advance\ypos by\halflength
       \advance\ypos by\gap
       \put(\xpos,\ypos){\evector(0,1){\halflength}}%
   \fi
\fi
}

\def\putmorphism(#1)(#2,#3)[#4`#5`#6]#7#8#9{{%
\run #2
\rise #3
\ifnum\rise=0
  \puthmorphism(#1)[#4`#5`#6]{#7}{#8}#9%
\else\ifnum\run=0
  \putvmorphism(#1)[#4`#5`#6]{#7}{#8}#9%
\else
\setpos(#1)%
\arrowlength #7
\arrowtype #8
\ifnum\run=0
\else\ifnum\rise=0
\else
\ifnum\run>0
    \coefa=1
\else
   \coefa=-1
\fi
\ifnum\arrowtype>0
   \coefb=0
   \coefc=-1
\else
   \coefb=\coefa
   \coefc=1
   \arrowtype=-\arrowtype
\fi
\width=2
\multiply \width by\run
\divide \width by\rise
\ifnum \width<0  \width=-\width\fi
\advance\width by60
\if l#9 \width=-\width\fi
\putbox(\xpos,\ypos){#4}%            %node 1
{\multiply \coefa by\arrowlength%      %node 2
\advance\xpos by\coefa
\multiply \coefa by\rise
\divide \coefa by\run
\advance \ypos by\coefa
\putbox(\xpos,\ypos){#5} }%
{\multiply \coefa by\arrowlength%      %label
\divide \coefa by2
\advance \xpos by\coefa
\advance \xpos by\width
\multiply \coefa by\rise
\divide \coefa by\run
\advance \ypos by\coefa
\if l#9%
   \putrbox(\xpos,\ypos){#6}%
\else\if r#9%
   \putlbox(\xpos,\ypos){#6}%
\fi\fi }%
{\multiply \rise by-\coefc%             %arrow
\multiply \run by-\coefc
\multiply \coefb by\arrowlength
\advance \xpos by\coefb
\multiply \coefb by\rise
\divide \coefb by\run
\advance \ypos by\coefb
\multiply \coefc by70
\advance \ypos by\coefc
\multiply \coefc by\run
\divide \coefc by\rise
\advance \xpos by\coefc
\multiply \coefa by140
\multiply \coefa by\run
\divide \coefa by\rise
\advance \arrowlength by\coefa
\ifcase\arrowtype
\or \put(\xpos,\ypos){\vector(\run,\rise){\arrowlength}}%
\or \put(\xpos,\ypos){\mvector(\run,\rise){\arrowlength}}%
\or \put(\xpos,\ypos){\evector(\run,\rise){\arrowlength}}%
\fi}\fi\fi\fi\fi}}

\newcount\numbdashes \newcount\lengthdash \newcount\increment

\def\howmanydashes{% Actually returns both number and length
\numbdashes=\arrowlength \lengthdash=40
\divide\numbdashes by \lengthdash
\lengthdash=\arrowlength
\divide\lengthdash by \numbdashes
\increment=\lengthdash
\multiply\lengthdash by 3
\divide\lengthdash by 5
}

\def\putdashvector(#1)(#2,#3)#4#5{%
\ifnum#3=0 \putdashhvector(#1){#4}#5
\else
\ifnum#2=0
\putdashvvector(#1){#4}#5\fi\fi}

\def\putdashhvector(#1,#2)#3#4{{%
\arrowlength=#3 \howmanydashes
\multiput(#1,#2)(\increment,0){\numbdashes}%
{\vrule height .4pt width \lengthdash\unitlength}
\arrowtype=#4 \xpos=#1
\ifnum\arrowtype<0 \advance\arrowtype by 7 \fi
\ifcase\arrowtype
\or \advance\xpos by 10
    \put(\xpos,#2){\vector(-1,0){\lengthdash}}
    \advance\xpos by 40
    \put(\xpos,#2){\vector(-1,0){\lengthdash}}
\or \advance \xpos by 10
    \put(\xpos,#2){\vector(-1,0){\lengthdash}}
    \advance\xpos by  \arrowlength
    \advance\xpos by  -50
    \put(\xpos,#2){\vector(-1,0){\lengthdash}}
\or \advance\xpos by 10
    \put(\xpos,#2){\vector(-1,0){\lengthdash}}
\or \advance\xpos by \arrowlength
    \advance\xpos by -\lengthdash
    \put(\xpos,#2){\vector(1,0){\lengthdash}}
\or {\advance\xpos by 10
    \put(\xpos,#2){\vector(1,0){\lengthdash}}}
    \advance\xpos by \arrowlength
    \advance\xpos by -\lengthdash
    \put(\xpos,#2){\vector(1,0){\lengthdash}}
\or \advance\xpos by \arrowlength
    \advance\xpos by -\lengthdash
    \put(\xpos,#2){\vector(1,0){\lengthdash}}
    \advance\xpos by -40
    \put(\xpos,#2){\vector(1,0){\lengthdash}}
   \fi
}}

\def\putdashvvector(#1,#2)#3#4{{%
\arrowlength=#3 \howmanydashes
\ypos=#2 \advance\ypos by -\arrowlength
\multiput(#1,#2)(0,\increment){\numbdashes}%
    {\vrule width .4pt height \lengthdash\unitlength}
\arrowtype=#4 \ypos=#2
\ifnum\arrowtype<0 \advance\arrowtype by 7 \fi
\ifcase\arrowtype
\or \advance\ypos by \arrowlength \advance\ypos by -40
    \put(#1,\ypos){\vector(0,1){\lengthdash}}
    \advance\ypos by -40
    \put(#1,\ypos){\vector(0,1){\lengthdash}}
\or \advance\ypos by 10
    \put(#1,\ypos){\vector(0,1){\lengthdash}}
    \advance\ypos by \arrowlength \advance\ypos by -40
    \put(#1,\ypos){\vector(0,1){\lengthdash}}
\or \advance\ypos by \arrowlength \advance\ypos by -40
    \put(#1,\ypos){\vector(0,1){\lengthdash}}
\or \advance\ypos by 10
    \put(#1,\ypos){\vector(0,-1){\lengthdash}}
\or \advance\ypos by 10
    \put(#1,\ypos){\vector(0,-1){\lengthdash}}
    \advance\ypos by \arrowlength \advance\ypos by -40
    \put(#1,\ypos){\vector(0,-1){\lengthdash}}
\or \advance\ypos by 10
    \put(#1,\ypos){\vector(0,-1){\lengthdash}}
    \advance\ypos by 40
    \put(#1,\ypos){\vector(0,-1){\lengthdash}}
\fi
}}

\def\puthmorphism(#1,#2)[#3`#4`#5]#6#7#8{{%
\xpos #1
\ypos #2
\width #6
\arrowlength #6
\arrowtype=#7
\putbox(\xpos,\ypos){#3\vphantom{#4}}%
{\advance \xpos by\arrowlength
\putbox(\xpos,\ypos){\vphantom{#3}#4}}%
\horsize{\tempcounta}{#3}%
\horsize{\tempcountb}{#4}%
\divide \tempcounta by2
\divide \tempcountb by2
\advance \tempcounta by30
\advance \tempcountb by30
\advance \xpos by\tempcounta
\advance \arrowlength by-\tempcounta
\advance \arrowlength by-\tempcountb
\putvector(\xpos,\ypos)(1,0)\arrowlength\arrowtype
\divide \arrowlength by2
\advance \xpos by\arrowlength
\vertsize{\tempcounta}{#5}%
\divide\tempcounta by2
\advance \tempcounta by20
\if a#8 %
   \advance \ypos by\tempcounta
   \putbox(\xpos,\ypos){#5}%
\else
   \advance \ypos by-\tempcounta
   \putbox(\xpos,\ypos){#5}%
\fi}}

\def\putvmorphism(#1,#2)[#3`#4`#5]#6#7#8{{%
\xpos #1
\ypos #2
\arrowlength #6
\arrowtype #7
\settowidth{\xlen}{$#5$}%
\putbox(\xpos,\ypos){#3}%
{\advance \ypos by-\arrowlength
\putbox(\xpos,\ypos){#4}}%
{\advance\arrowlength by-140
\advance \ypos by-70
\ifdim\xlen>0pt
   \if m#8%
      \putsplitvector(\xpos,\ypos)\arrowlength\arrowtype
   \else
   \putvector(\xpos,\ypos)(0,-1)\arrowlength\arrowtype
   \fi
\else
   \putvector(\xpos,\ypos)(0,-1)\arrowlength\arrowtype
\fi}%
\ifdim\xlen>0pt
   \divide \arrowlength by2
   \advance\ypos by-\arrowlength
   \if l#8%
      \advance \xpos by-40
      \putrbox(\xpos,\ypos){#5}%
   \else\if r#8%
      \advance \xpos by40
      \putlbox(\xpos,\ypos){#5}%
   \else
      \putbox(\xpos,\ypos){#5}%
   \fi\fi
\fi
}}

\def\putsquarep<#1>(#2)[#3;#4`#5`#6`#7]{{%
\setsqparms[#1]%
\setpos(#2)%
\settokens`#3`%
\puthmorphism(\xpos,\ypos)[\tokenc`\tokend`{#7}]{\width}{\arrowtyped}b%
\advance\ypos by \height
\puthmorphism(\xpos,\ypos)[\tokena`\tokenb`{#4}]{\width}{\arrowtypea}a%
\putvmorphism(\xpos,\ypos)[``{#5}]{\height}{\arrowtypeb}l%
\advance\xpos by \width
\putvmorphism(\xpos,\ypos)[``{#6}]{\height}{\arrowtypec}r%
}}

\def\putsquare{\@ifnextchar <{\putsquarep}{\putsquarep%
   <\arrowtypea`\arrowtypeb`\arrowtypec`\arrowtyped;\width`\height>}}
\def\square{\@ifnextchar< {\squarep}{\squarep
   <\arrowtypea`\arrowtypeb`\arrowtypec`\arrowtyped;\width`\height>}}
\def\squarep<#1>[#2`#3`#4`#5;#6`#7`#8`#9]{{%       %     #2------>#3
\setsqparms[#1]%                                   %      |       |
\diagram%                                          %      |       |
\putsquarep<\arrowtypea`\arrowtypeb`\arrowtypec`%  %    #7|       |#8
\arrowtyped;\width`\height>%                       %      |       |
(0,0)[#2`#3`#4`{#5};#6`#7`#8`{#9}]%                %      |       |
\enddiagram%                                       %      v       v
}}                                                 %     #4------>#5
\def\putptrianglep<#1>(#2,#3)[#4`#5`#6;#7`#8`#9]{{%
\settriparms[#1]%
\xpos=#2 \ypos=#3
\advance\ypos by \height
\puthmorphism(\xpos,\ypos)[#4`#5`{#7}]{\height}{\arrowtypea}a%
\putvmorphism(\xpos,\ypos)[`#6`{#8}]{\height}{\arrowtypeb}l%
\advance\xpos by\height
\putmorphism(\xpos,\ypos)(-1,-1)[``{#9}]{\height}{\arrowtypec}r%
}}

\def\putptriangle{\@ifnextchar <{\putptrianglep}{\putptrianglep
   <\arrowtypea`\arrowtypeb`\arrowtypec;\height>}}
\def\ptriangle{\@ifnextchar <{\ptrianglep}{\ptrianglep
   <\arrowtypea`\arrowtypeb`\arrowtypec;\height>}}
\def\ptrianglep<#1>[#2`#3`#4;#5`#6`#7]{{%%    %      #2----->#3
\settriparms[#1]%                             %      |      /
\diagram%                                     %      |     /
\putptrianglep<\arrowtypea`\arrowtypeb`%      %    #6|    /#7
\arrowtypec;\height>%                         %      |   /
(0,0)[#2`#3`#4;#5`#6`{#7}]%                   %      |  /
\enddiagram%%                                 %      v v
}}                                            %      #4

\def\putqtrianglep<#1>(#2,#3)[#4`#5`#6;#7`#8`#9]{{%
\settriparms[#1]%
\xpos=#2 \ypos=#3
\advance\ypos by\height
\puthmorphism(\xpos,\ypos)[#4`#5`{#7}]{\height}{\arrowtypea}a%
\putmorphism(\xpos,\ypos)(1,-1)[``{#8}]{\height}{\arrowtypeb}l%
\advance\xpos by\height
\putvmorphism(\xpos,\ypos)[`#6`{#9}]{\height}{\arrowtypec}r%
}}

\def\putqtriangle{\@ifnextchar <{\putqtrianglep}{\putqtrianglep
   <\arrowtypea`\arrowtypeb`\arrowtypec;\height>}}
\def\qtriangle{\@ifnextchar <{\qtrianglep}{\qtrianglep
   <\arrowtypea`\arrowtypeb`\arrowtypec;\height>}}
\def\qtrianglep<#1>[#2`#3`#4;#5`#6`#7]{{%%    %        #2----->#3
\settriparms[#1]%                             %         \      |
\width=\height                                %          \     |
\diagram%                                     %         #6\    |#7
\putqtrianglep<\arrowtypea`\arrowtypeb`%      %            \   |
\arrowtypec;\height>%                         %             \  |
(0,0)[#2`#3`#4;#5`#6`{#7}]%                   %              v v
\enddiagram%%                                 %               #4
}}

\def\putdtrianglep<#1>(#2,#3)[#4`#5`#6;#7`#8`#9]{{%
\settriparms[#1]%
\xpos=#2 \ypos=#3
\puthmorphism(\xpos,\ypos)[#5`#6`{#9}]{\height}{\arrowtypec}b%
\advance\xpos by \height \advance\ypos by\height
\putmorphism(\xpos,\ypos)(-1,-1)[``{#7}]{\height}{\arrowtypea}l%
\putvmorphism(\xpos,\ypos)[#4``{#8}]{\height}{\arrowtypeb}r%
}}

\def\putdtriangle{\@ifnextchar <{\putdtrianglep}{\putdtrianglep
   <\arrowtypea`\arrowtypeb`\arrowtypec;\height>}}
\def\dtriangle{\@ifnextchar <{\dtrianglep}{\dtrianglep
   <\arrowtypea`\arrowtypeb`\arrowtypec;\height>}}
\def\dtrianglep<#1>[#2`#3`#4;#5`#6`#7]{{%%    %                  / |
\settriparms[#1]%                             %                 /  |
\width=\height                                %              #5/   |#6
\diagram%                                     %               /    |
\putdtrianglep<\arrowtypea`\arrowtypeb`%      %              /     |
\arrowtypec;\height>%                         %             v      v
(0,0)[#2`#3`#4;#5`#6`{#7}]%                   %            #3----->#4
\enddiagram%%                                 %                #7
}}

\def\putbtrianglep<#1>(#2,#3)[#4`#5`#6;#7`#8`#9]{{%
\settriparms[#1]%
\xpos=#2 \ypos=#3
\puthmorphism(\xpos,\ypos)[#5`#6`{#9}]{\height}{\arrowtypec}b%
\advance\ypos by\height
\putmorphism(\xpos,\ypos)(1,-1)[``{#8}]{\height}{\arrowtypeb}r%
\putvmorphism(\xpos,\ypos)[#4``{#7}]{\height}{\arrowtypea}l%
}}

\def\putbtriangle{\@ifnextchar <{\putbtrianglep}{\putbtrianglep
   <\arrowtypea`\arrowtypeb`\arrowtypec;\height>}}
\def\btriangle{\@ifnextchar <{\btrianglep}{\btrianglep
   <\arrowtypea`\arrowtypeb`\arrowtypec;\height>}}
\def\btrianglep<#1>[#2`#3`#4;#5`#6`#7]{{%%   %              | \
\settriparms[#1]%                            %              |  \
\width=\height                               %            #5|   \#6
\diagram%                                    %              |    \
\putbtrianglep<\arrowtypea`\arrowtypeb`%     %              |     \
\arrowtypec;\height>%                        %              v      v
(0,0)[#2`#3`#4;#5`#6`{#7}]%                  %              #3----->#4
\enddiagram%%                                %                 #7
}}

\def\putAtrianglep<#1>(#2,#3)[#4`#5`#6;#7`#8`#9]{{%
\settriparms[#1]%
\xpos=#2 \ypos=#3
{\multiply \height by2
\puthmorphism(\xpos,\ypos)[#5`#6`{#9}]{\height}{\arrowtypec}b}%
\advance\xpos by\height \advance\ypos by\height
\putmorphism(\xpos,\ypos)(-1,-1)[#4``{#7}]{\height}{\arrowtypea}l%
\putmorphism(\xpos,\ypos)(1,-1)[``{#8}]{\height}{\arrowtypeb}r%
}}

\def\putAtriangle{\@ifnextchar <{\putAtrianglep}{\putAtrianglep
   <\arrowtypea`\arrowtypeb`\arrowtypec;\height>}}
\def\Atriangle{\@ifnextchar <{\Atrianglep}{\Atrianglep
   <\arrowtypea`\arrowtypeb`\arrowtypec;\height>}}
\def\Atrianglep<#1>[#2`#3`#4;#5`#6`#7]{{%%         %         /   \
\settriparms[#1]%                                  %        /     \
\width=\height                                     %     #5/       \#6
\diagram%                                          %      /         \
\putAtrianglep<\arrowtypea`\arrowtypeb`%           %     /           \
\arrowtypec;\height>%                              %    v             v
(0,0)[#2`#3`#4;#5`#6`{#7}]%                        %   #3------------>#4
\enddiagram%%                                      %          #7
}}

\def\putAtrianglepairp<#1>(#2)[#3;#4`#5`#6`#7`#8]{{%
\settripairparms[#1]%
\setpos(#2)%
\settokens`#3`%
\puthmorphism(\xpos,\ypos)[\tokenb`\tokenc`{#7}]{\height}{\arrowtyped}b%
\advance\xpos by\height
\puthmorphism(\xpos,\ypos)[\phantom{\tokenc}`\tokend`{#8}]%
{\height}{\arrowtypee}b%
\advance\ypos by\height
\putmorphism(\xpos,\ypos)(-1,-1)[\tokena``{#4}]{\height}{\arrowtypea}l%
\putvmorphism(\xpos,\ypos)[``{#5}]{\height}{\arrowtypeb}m%
\putmorphism(\xpos,\ypos)(1,-1)[``{#6}]{\height}{\arrowtypec}r%
}}

\def\putAtrianglepair{\@ifnextchar <{\putAtrianglepairp}{\putAtrianglepairp%
   <\arrowtypea`\arrowtypeb`\arrowtypec`\arrowtyped`\arrowtypee;\height>}}
\def\Atrianglepair{\@ifnextchar <{\Atrianglepairp}{\Atrianglepairp%
   <\arrowtypea`\arrowtypeb`\arrowtypec`\arrowtyped`\arrowtypee;\height>}}

\def\Atrianglepairp<#1>[#2;#3`#4`#5`#6`#7]{{%           %  #2a
\settripairparms[#1]%                         %           / | \
\settokens`#2`%                               %          /  |  \
\width=\height                                %       #3/  #4   \#5
\diagram%                                     %        /    |    \
\putAtrianglepairp                            %       /     |     \
<\arrowtypea`\arrowtypeb`\arrowtypec`%        %      v      v      v
\arrowtyped`\arrowtypee;\height>%             %     #2b---->#2c---->#2d
(0,0)[{#2};#3`#4`#5`#6`{#7}]%                 %         #6     #7
\enddiagram%%
}}

\def\putVtrianglep<#1>(#2,#3)[#4`#5`#6;#7`#8`#9]{{%
\settriparms[#1]%
\xpos=#2 \ypos=#3
\advance\ypos by\height
{\multiply\height by2
\puthmorphism(\xpos,\ypos)[#4`#5`{#7}]{\height}{\arrowtypea}a}%
\putmorphism(\xpos,\ypos)(1,-1)[`#6`{#8}]{\height}{\arrowtypeb}l%
\advance\xpos by\height
\advance\xpos by\height
\putmorphism(\xpos,\ypos)(-1,-1)[``{#9}]{\height}{\arrowtypec}r%
}}

\def\putVtriangle{\@ifnextchar <{\putVtrianglep}{\putVtrianglep
   <\arrowtypea`\arrowtypeb`\arrowtypec;\height>}}
\def\Vtriangle{\@ifnextchar <{\Vtrianglep}{\Vtrianglep
   <\arrowtypea`\arrowtypeb`\arrowtypec;\height>}}
\def\Vtrianglep<#1>[#2`#3`#4;#5`#6`#7]{{%%     %        #2------------->#3
\settriparms[#1]%                              %         \             /
\width=\height                                 %          \           /
\diagram%                                      %         #6\         /#7
\putVtrianglep<\arrowtypea`\arrowtypeb`%       %            \       /
\arrowtypec;\height>%                          %             \     /
(0,0)[#2`#3`#4;#5`#6`{#7}]%                    %              v   v
\enddiagram%%                                  %               #4
}}

\def\putVtrianglepairp<#1>(#2)[#3;#4`#5`#6`#7`#8]{{
\settripairparms[#1]%
\setpos(#2)%
\settokens`#3`%
\advance\ypos by\height
\putmorphism(\xpos,\ypos)(1,-1)[`\tokend`{#6}]{\height}{\arrowtypec}l%
\puthmorphism(\xpos,\ypos)[\tokena`\tokenb`{#4}]{\height}{\arrowtypea}a%
\advance\xpos by\height
\puthmorphism(\xpos,\ypos)[\phantom{\tokenb}`\tokenc`{#5}]%
{\height}{\arrowtypeb}a%
\putvmorphism(\xpos,\ypos)[``{#7}]{\height}{\arrowtyped}m%
\advance\xpos by\height
\putmorphism(\xpos,\ypos)(-1,-1)[``{#8}]{\height}{\arrowtypee}r%
}}

\def\putVtrianglepair{\@ifnextchar <{\putVtrianglepairp}{\putVtrianglepairp%
    <\arrowtypea`\arrowtypeb`\arrowtypec`\arrowtyped`\arrowtypee;\height>}}
\def\Vtrianglepair{\@ifnextchar <{\Vtrianglepairp}{\Vtrianglepairp%
    <\arrowtypea`\arrowtypeb`\arrowtypec`\arrowtyped`\arrowtypee;\height>}}
\def\Vtrianglepairp<#1>[#2;#3`#4`#5`#6`#7]{{%  %  #2a---->#2b---->#2c
\settripairparms[#1]%                          %   \      |      /
\settokens`#2`%                                %    \     |     /
\diagram%                                      %   #5\   #6    /#7
\putVtrianglepairp                             %      \   |   /
<\arrowtypea`\arrowtypeb`\arrowtypec`%         %       \  |  /
\arrowtyped`\arrowtypee;\height>%              %        v v v
(0,0)[{#2};#3`#4`#5`#6`{#7}]%                  %         #2d
\enddiagram%%
}}

\def\putCtrianglep<#1>(#2,#3)[#4`#5`#6;#7`#8`#9]{{%
\settriparms[#1]%
\xpos=#2 \ypos=#3
\advance\ypos by\height
\putmorphism(\xpos,\ypos)(1,-1)[``{#9}]{\height}{\arrowtypec}l%
\advance\xpos by\height
\advance\ypos by\height
\putmorphism(\xpos,\ypos)(-1,-1)[#4`#5`{#7}]{\height}{\arrowtypea}l%
{\multiply\height by 2
\putvmorphism(\xpos,\ypos)[`#6`{#8}]{\height}{\arrowtypeb}r}%
}}

\def\putCtriangle{\@ifnextchar <{\putCtrianglep}{\putCtrianglep
    <\arrowtypea`\arrowtypeb`\arrowtypec;\height>}}
\def\Ctriangle{\@ifnextchar <{\Ctrianglep}{\Ctrianglep
    <\arrowtypea`\arrowtypeb`\arrowtypec;\height>}}
\def\Ctrianglep<#1>[#2`#3`#4;#5`#6`#7]{{%%   %                / |
\settriparms[#1]%                            %             #5/  |
\width=\height                               %              /   |
\diagram%                                    %             v    |
\putCtrianglep<\arrowtypea`\arrowtypeb`%     %           #3     |#6
\arrowtypec;\height>%                        %             \    |
(0,0)[#2`#3`#4;#5`#6`{#7}]%                  %            #7\   |
\enddiagram%%                                %               \  |
}}                                           %                v v
\def\putDtrianglep<#1>(#2,#3)[#4`#5`#6;#7`#8`#9]{{%
\settriparms[#1]%
\xpos=#2 \ypos=#3
\advance\xpos by\height \advance\ypos by\height
\putmorphism(\xpos,\ypos)(-1,-1)[``{#9}]{\height}{\arrowtypec}r%
\advance\xpos by-\height \advance\ypos by\height
\putmorphism(\xpos,\ypos)(1,-1)[`#5`{#8}]{\height}{\arrowtypeb}r%
{\multiply\height by 2
\putvmorphism(\xpos,\ypos)[#4`#6`{#7}]{\height}{\arrowtypea}l}%
}}

\def\putDtriangle{\@ifnextchar <{\putDtrianglep}{\putDtrianglep
    <\arrowtypea`\arrowtypeb`\arrowtypec;\height>}}
\def\Dtriangle{\@ifnextchar <{\Dtrianglep}{\Dtrianglep
   <\arrowtypea`\arrowtypeb`\arrowtypec;\height>}}
\def\Dtrianglep<#1>[#2`#3`#4;#5`#6`#7]{{%%  %          | \
\settriparms[#1]%                           %          |  \#6
\width=\height                              %          |   \
\diagram%                                   %          |    v
\putDtrianglep<\arrowtypea`\arrowtypeb`%    %        #5|    #3
\arrowtypec;\height>%                       %          |    /
(0,0)[#2`#3`#4;#5`#6`{#7}]%                 %          |   /#7
\enddiagram%%                               %          |  /
}}                                          %          v v
\def\setrecparms[#1`#2]{\width=#1 \height=#2}%
\def\recursep<#1`#2>[#3;#4`#5`#6`#7`#8]{{\m@th
\width=#1 \height=#2
\settokens`#3`
\settowidth{\tempdimen}{$\tokena$}
\ifdim\tempdimen=0pt
  \savebox{\tempboxa}{\hbox{$\tokenb$}}%
  \savebox{\tempboxb}{\hbox{$\tokend$}}%
  \savebox{\tempboxc}{\hbox{$#6$}}%
\else
  \savebox{\tempboxa}{\hbox{$\hbox{$\tokena$}\times\hbox{$\tokenb$}$}}%
  \savebox{\tempboxb}{\hbox{$\hbox{$\tokena$}\times\hbox{$\tokend$}$}}%
  \savebox{\tempboxc}{\hbox{$\hbox{$\tokena$}\times\hbox{$#6$}$}}%
\fi
\ypos=\height
\divide\ypos by 2
\xpos=\ypos
\advance\xpos by \width
\bfig
\putCtrianglep<-1`1`1;\ypos>(0,0)[`\tokenc`;#5`#6`{#7}]%
\puthmorphism(\ypos,0)[\tokend`\usebox{\tempboxb}`{#8}]{\width}{-1}b%
\puthmorphism(\ypos,\height)[\tokenb`\usebox{\tempboxa}`{#4}]{\width}{-1}a%
\advance\ypos by \width
\putvmorphism(\ypos,\height)[``\usebox{\tempboxc}]{\height}1r%
\efig
}}

\def\recurse{\@ifnextchar <{\recursep}{\recursep<\width`\height>}}

\def\puttwohmorphisms(#1,#2)[#3`#4;#5`#6]#7#8#9{{%
\puthmorphism(#1,#2)[#3`#4`]{#7}0a
\ypos=#2
\advance\ypos by 20
\puthmorphism(#1,\ypos)[\phantom{#3}`\phantom{#4}`#5]{#7}{#8}a
\advance\ypos by -40
\puthmorphism(#1,\ypos)[\phantom{#3}`\phantom{#4}`#6]{#7}{#9}b
}}

\def\puttwovmorphisms(#1,#2)[#3`#4;#5`#6]#7#8#9{{%
\putvmorphism(#1,#2)[#3`#4`]{#7}0a
\xpos=#1
\advance\xpos by -20
\putvmorphism(\xpos,#2)[\phantom{#3}`\phantom{#4}`#5]{#7}{#8}l
\advance\xpos by 40
\putvmorphism(\xpos,#2)[\phantom{#3}`\phantom{#4}`#6]{#7}{#9}r
}}

\def\puthcoequalizer(#1)[#2`#3`#4;#5`#6`#7]#8#9{{%
\setpos(#1)%
\puttwohmorphisms(\xpos,\ypos)[#2`#3;#5`#6]{#8}11%
\advance\xpos by #8
\puthmorphism(\xpos,\ypos)[\phantom{#3}`#4`#7]{#8}1{#9}
}}

\def\putvcoequalizer(#1)[#2`#3`#4;#5`#6`#7]#8#9{{%
\setpos(#1)%
\puttwovmorphisms(\xpos,\ypos)[#2`#3;#5`#6]{#8}11%
\advance\ypos by -#8
\putvmorphism(\xpos,\ypos)[\phantom{#3}`#4`#7]{#8}1{#9}
}}

\def\putthreehmorphisms(#1)[#2`#3;#4`#5`#6]#7(#8)#9{{%
\setpos(#1) \settypes(#8)
\if a#9 %
     \vertsize{\tempcounta}{#5}%
     \vertsize{\tempcountb}{#6}%
     \ifnum \tempcounta<\tempcountb \tempcounta=\tempcountb \fi
\else
     \vertsize{\tempcounta}{#4}%
     \vertsize{\tempcountb}{#5}%
     \ifnum \tempcounta<\tempcountb \tempcounta=\tempcountb \fi
\fi
\advance \tempcounta by 60
\puthmorphism(\xpos,\ypos)[#2`#3`#5]{#7}{\arrowtypeb}{#9}
\advance\ypos by \tempcounta
\puthmorphism(\xpos,\ypos)[\phantom{#2}`\phantom{#3}`#4]{#7}{\arrowtypea}{#9}
\advance\ypos by -\tempcounta \advance\ypos by -\tempcounta
\puthmorphism(\xpos,\ypos)[\phantom{#2}`\phantom{#3}`#6]{#7}{\arrowtypec}{#9}
}}

\def\setarrowtoks[#1`#2`#3`#4`#5`#6]{%
\def\toka{#1}
\def\tokb{#2}
\def\tokc{#3}
\def\tokd{#4}
\def\toke{#5}
\def\tokf{#6}
}
\def\hex{\@ifnextchar <{\hexp}{\hexp<1000`400>}}
\def\hexp<#1`#2>[#3`#4`#5`#6`#7`#8;#9]{%
\setarrowtoks[#9]
\yext=#2 \advance \yext by #2
\xext=#1 \advance\xext by \yext
\bfig
\putCtriangle<-1`0`1;#2>(0,0)[`#5`;\tokb``\tokd]
\xext=#1 \yext=#2 \advance \yext by #2
\putsquare<1`0`0`1;\xext`\yext>(#2,0)[#3`#4`#7`#8;\toka```\tokf]
\advance \xext by #2
\putDtriangle<0`1`-1;#2>(\xext,0)[`#6`;`\tokc`\toke]
\efig
}
\makeatother